\documentstyle[11pt]{article}
\newtheorem {th}{Theorem}[section]
\newtheorem {lem}[th]{Lemma}
\newtheorem {pr}[th]{Proposition}
\newtheorem {cor}[th]{Corollary}

\def\wt{\widetilde}
\def\wqt{\wt{q}}
\def\Cox{\hfill \Box}
\def\R{{\bf R}} 
\def\wt{\widetilde} 
\def\grad{\bigtriangledown}
\def\lap{\bigtriangledown^2}
 
\def\vv{{\bf v}}
\def\pp{{\bf p}}
\def\uu{{\bf u}}
\def\ww{{\bf w}}
\def\xx{{\bf x}}
\def\yy{{\bf y}}
\def\zz{{\bf z}}

\def\deq{:=}

\def\dd{\delta}
\def\ee{\epsilon}
\def\E{{\bf{E}}}
\def\P{{\bf{P}}}
\def\N{\hbox{I\kern-.2em\hbox{N}}}
\def\R{\hbox{I\kern-.2em\hbox{R}}}

\def\Z{{\bf{Z}}}
\def\B{{\cal{B}}}
\def\F{{\cal{F}}}
\def\C{{\cal {C}}}

\def\cR{{\cal {R}}}

\def\cl{c_{\rm cyl}}
\def\clo{\cl'}
\def\cll{{\tilde{c}}}

\def\cllo{c_{\rm sum}}
\def\clll{c_q}
\def\clool{{c_{\rm density}}}
\def\clolo{{c_{\rm density}'}}
\def\cloll{{c_{\rm density}''}}
\def\chbounds{c_*}
\def\clloo{{c_1}}
\def\cllol{{c_{\rm upper}}}
\def\clllo{{c_{\rm upper}'}}
\def\cllll{{c_{\rm upper}''}}
\def\cloooo{{c_{\rm upper}^*}}
\def\cloool{\zeta}
\def\cloolo{M}

\def\|{\, | \, }
\def\one{{\bf 1}}

\def\K{{\cal{K}}}
\def\W{{\cal{W}}}

\begin{document}
\begin{titlepage}
\begin{center}
{\large \bf SETS AVOIDED BY BROWNIAN MOTION} \\
\end{center}
\vspace{5ex}
\begin{center} {\bf
Omer Adelman\footnote{Universit\'e Pierre et Marie Curie, 
Laboratoire de Probabilit\'es, 4, Place Jussieu -- Tour 56, 
75252 Paris Cedex 05, FRANCE} \\
Krzysztof Burdzy\footnote{University of Washington, 
Department of Mathematics, Box 354350, Seattle, WA 98195--4350}$^,$
\footnote{Research supported in part by National Science Foundation Grant 
\# DMS 9322689} \\
Robin Pemantle \footnote{Department of Mathematics, University of 
Wisconsin-Madison, Van Vleck Hall, 480 Lincoln Drive, Madison, WI 53706}$^,$
\footnote{Research supported in part by National Science Foundation Grant 
\# DMS 9300191, by a Sloan Foundation Fellowship and by a Presidential
Faculty Fellowship} ~\\
}\end{center}

\vfill

ABSTRACT:

A fixed 2-dimensional projection of a 3-dimensional Brownian 
motion is almost surely neighborhood recurrent; is this simultaneously 
true of all the 2-dimensional projections with probability one?
Equivalently: 3-dimensional Brownian motion hits any infinite 
cylinder with probability one; does it hit all cylinders?  This
papers shows that the answer is no.  Brownian motion in three dimensions
avoids random cylinders and in fact avoids bodies of revolution that
grow almost as fast as cones.  

\vfill



\end{titlepage}

\section{Introduction}

\setcounter{equation}{0}

Let $\{ B(t) : 0 \leq t < \infty \}$ be a Brownian motion started
 from the origin in three dimensions, with coordinates $(X(t), Y(t),
Z(t))$ defined on $(\Omega, \F (t) , \P )$; we will use $\P_\vv$
to denote the law of $B(t)$ translated by $\vv$.  Any projection of
$\{ B(t) \}$ onto a plane is a version of a 2-dimensional Brownian
motion and its almost sure properties are well known: it is
neighborhood recurrent, its range is two-dimensional with exact 
Hausdorff gauge $x^2 \log (1/x) \log \log \log (1/x)$; the list goes on.  
Some of these properties are known to hold uniformly over all
projections, while others fail (necessarily on a set of projections
of measure zero, by Fubini's theorem).  What about neighborhood
recurrence: is this a property inherited simultaneously by all
projections of $\{ B(t) \}$?  An equivalent question is: 
\begin{center}
Does
$\{ B(t) \}$ with probability one intersect every infinite cylinder?
\end{center}

In this paper we give a negative answer: with probability one, there
are random cylinders disjoint from the range of a 3-dimensional
Brownian motion.  In fact we show more.  Let $f$ be a strictly positive
increasing function on $\R^+$ and let $\C_f$ be the set or {\em thorn}
$$\left \{ (x,y,z) \in \R^3 : x^2 + y^2 + z^2 \geq 1 \mbox{ and }
   \sqrt{x^2 + y^2} \leq f(|z|) \right \} .$$  
Say that Brownian motion {\em avoids $f$-thorns} if there is with 
probability one a random set congruent to $\C_f$ avoided by Brownian
motion.  A zero-one law holds, so the alternative is that with probability
one Brownian motion intersects all sets congruent to $\C_f$.
Our main results are contained in Theorems~\ref{th converse} 
and~\ref{th integral test} below: under an integral condition on $f$, 
satisfied for example when $f(z) = z / \exp (\log^{1/2 + \ee} z)$, 
Brownian motion avoids $f$-thorns; moreover, in this case the set
of directions of axes of $f$-thorns avoided by Brownian motion has
Hausdorff dimension 2, with positive probability. On the other hand, if 
$f(z) = z / \exp (c \log^{1/2} z)$ for sufficiently small $c$, 
then Brownian motion does not avoid $f$-thorns, a.s.

{\em Remarks:} 

\noindent{1.}  It would have been equally natural to consider
one-sided thorns, $\C_f \cap \{ (x,y,z) : z \geq 0 \}$, but
there seems to be little difference since we cannot find an
$f$ for which Brownian motion intersects all two-sided
$f$-thorns but misses some one-sided $f$-thorns.

\noindent{2.}  One original motivation for this question was to shed
some light on the complement of the Wiener sausage 
$\W \deq \{ B(t) + \xx : t \in \R^+ , |\xx| \leq 1\}$.  For example,
we do not know an elementary proof that $\R^3 \setminus \W$ has
an unbounded connected component.  This follows from the weakest 
of our avoidance results.  An elementary argument, based on the
existence of arbtrarily large values of $t$ for which
$\sup_{s \leq t} |B(s)| - \inf_{s \geq t} |B(s)|$ is
smaller than any arbitrary fixed positive number
(see~\cite{AdSh} or~\cite[Proposition~1]{Bu}), shows that there
must be at most one unbounded component.  

\noindent{3. } The notion of properties holding uniformly over planar
projections of higher dimensional Brownian motion is similar to the notion
of {\em quasi-everywhere} properties of the Brownian path, that is,
properties that w.~p.~1 hold simultaneously for every cross section of the
Brownian Sheet.  See for example~\cite{Fu} or~\cite{Pen}.

We now briefly outline the arguments, setting forth notation that
will be used throughout.
\begin{quote}
{\bf Notation} For any unit vector $\vv\in\R^3$, let $\C_\vv = \C_{f , \vv}$ 
denote the image of $\C_f$ under any 
origin-preserving rotation mapping $(0,0,1)$ to $\vv$.
Usually $f$ will be fixed and will be dropped from the notation. 
Let $\vv_\theta$ denote $(\sin (\theta) , 0 , \cos (\theta))$ and
let $\C_\theta$ denote $\C_{\vv_\theta}$.  For any set $A$, let 
$\tau_A$ denote the time Brownian motion first hits the set $A$.
Let $\B (\xx , L)$ denote the ball of radius $L$ about the point $\xx$, 
let $\B_L$ denote $\B (0 , L)$, and let $\tau_L$ be shorthand for
$\tau_{\partial \B_L}$.   Let 
$$q(L) = \P (\tau_L < \tau_\C)$$ 
be the probability that Brownian motion reaches modulus $L$ before
hitting the $f$-thorn.  Let 
$$q(L , \theta) = \P (\tau_L < \tau_\C \wedge \tau_{\C_\theta})$$
be the probability that Brownian motion reaches modulus $L$ before
hitting either of two $f$-thorns separated by an angle of $\theta$.  
Write $\mu_L$ for the hitting subprobability measure on $\partial \B_L$
of Brownian motion absorbed by $\C$, so that for 
$A \subseteq \partial \B_L$, 
$\mu_L (A) = \P (\tau_L < \tau_\C , B(\tau_L) \in A)$.  Let
$\mu_{L , \theta}$ be the same for $\C \cup \C_\theta$:
$$\mu_{L , \theta} (A) = \P (\tau_L < \tau_\C \wedge \tau_{\C_\theta} ,
   B(\tau_L) \in A) .$$
\end{quote}

Theorem~\ref{th integral test} is proved by the second moment method and
the easier Theorem~\ref{th converse} is proved by a first moment 
estimate.  We first restrict our attention from all sets congruent to
$\C$ to only the rotations, $\C_\vv$: 
we will show that Brownian motion avoids
$f$-thorns if and only if with positive probability there is a
set $\C_\vv$ avoided by Brownian motion.  Let
$W_L$ be the measure of the set of all vectors $\vv$ in the unit sphere 
for which $\tau_L < \tau_{\C_\vv}$.  Estimates on $q(L)$ and
$q(L , \theta)$ yield estimates on $\E W_L$ and $\E W_L^2$.
When $\E W_L^2 / (\E W_L)^2$ is bounded, it follows that
$\liminf \P (W_L > 0) > 0$ and hence that Brownian motion
avoids $f$-thorns; when $\E W_L = o(f(z) / z)$, it follows that
$\P (W_L > 0) \rightarrow 0$ and hence that Brownian motion 
does not avoid $f$-thorns.  

All the work is in obtaining the estimates on $q(L)$ and, particularly, 
$q(L , \theta)$.  The remainder of the paper is organized as follows. 
The next section contains precise statements of the main results and
contains rigorous versions of the arguments mentioned above 
(zero-one laws, the first and second moment methods).  It also 
contains a proof of Theorem~\ref{th converse}, which requires very
little computation.  Section~3
contains proofs of the estimates on $q(L)$ and $q(L , \theta)$ in
the special case where $f(z) = z^\alpha$.  The reason for separating
this from the general case is that in the $z^\alpha$ case we have
reasonably accurate estimates of both $q(L)$ and $q(L , \theta)$.
While the boundedness of $\E W_L^2 / (\E W_L)^2$ in this case
is subsumed by our later results, these do not contain separate
estimates for $q(L)$ and $q(L , \theta)$, and we suspect that
the estimate on $q(L)$, Lemma~\ref{lem P(L)}, will be useful
in other contexts.  Section~4 begins the proof of 
Theorem~\ref{th integral test},
breaking it down into a series of lemmas.  Section~5 proves those lemmas
with soft proofs, Section~6 proves those lemmas involving manipulation
of Green's functions, and Section~7 gives the proofs that require geometric
analysis.  The most important ingredient in these last four sections
is the integration by parts device, Theorem~\ref{th clean green}, which
allows the computation of $U(L , \theta) := q(L , \theta) / q(L)^2$
without exact or even asymptotic knowledge of $q(L)$.  In addition 
to sharpening the dividing line between thorns that are avoided and 
thorns that are not, Theorem~\ref{th clean green} should be
useful in any situation where one wishes to estimate the probability
of simultaneously avoiding two sets.  

We will use many notions and results from the classical potential
theory and their probabilistic counterparts. A good presentation
of different aspects of this theory may be found in~\cite{Bass,KS,Port}.

\section{Main results}

\setcounter{equation}{0}

Write $\cR$ for the range of the Brownian motion $\{ B(t) \}$.  Throughout,
we let $g(z)$ denote the function $z / f(z)$.  From the fact that the radial
projection of Brownian motion onto the unit sphere is dense, we get the 
well known fact that Brownian motion cannot avoid cones:
\begin{th} \label{th 2.0}
If $f (z) = c z$ for some $c > 0$ then Brownian motion does not avoid
$f$-thorns.
\end{th}
$\Cox$

In the next section, we prove a first result in the other direction:
\begin{th} \label{th 2.3}
If $f(z) = z^\alpha$ for some $\alpha \in [0,1)$, then Brownian
motion avoids $f$-thorns.
\end{th}
In particular, when $\alpha = 0$ we recover the result first mentioned in the
introduction: some planar projections of 3-dimensional Brownian motion are
not neighborhood recurrent.  

Our sharpest non-avoidance result is:
\begin{th} \label{th converse}
If $f (z) = z / \exp (c \log^{1/2} z)$ for $c > 0$ sufficiently small, then 
Brownian motion does not avoid $f$-thorns.
\end{th}
Let
$ A = \{ \vv : |\vv| = 1 , \cR \cap \C_\vv = \emptyset \}$
be the set of directions of $f$-thorns
avoided by Brownian motion.
Our sharpest avoidance result is:
\begin{th} \label{th integral test}
Assume the following hypotheses on $f$ and on $g(z) := z / f(z)$:
\begin{eqnarray}
&& f(z) \mbox{ and } g(z) \mbox{ are increasing and tend to infinity as } 
z \rightarrow \infty, \label{eq hyp1} \\
&& g(1) \geq 2, \label{eq hyp2} \\
&& \mbox{the circle } A \mbox{ lies inside the region } |x| \leq f(z) ,
   \label{eq hyp3} 
\end{eqnarray}
where $A$ is the circle in the $z$-$x$ plane centered on the $z$ axis and
tangent to the graph $|x| = f(z)$ at the points $(z , \pm f(z))$.  If 
\begin{equation} \label{eq integral test}
\int_1^\infty {1 \over z \log^2 g(z)} \, dz < \infty ,
\end{equation} 
then Brownian motion avoids $f$-thorns, and in fact the set 
$ A = \{ \vv : |\vv| = 1 , \cR \cap \C_\vv = \emptyset \}$ of directions
of axes of $f$-thorns avoided by Brownian motion has Hasdorff dimension 2,
with positive probability.
\end{th}
\medskip
Note that the set 
$ A = \{ \vv : |\vv| = 1 , \cR \cap \C_\vv = \emptyset \}$ of directions
of axes of $f$-thorns avoided by Brownian motion
can be empty, with positive probability, for every non-trivial $f$.
\medskip

\noindent{\em Remarks on the hypotheses:}   Whether Brownian
motion avoids $f$-thorns is a monotone function of $f$ and does
not depend on the values of $f$ on any bounded interval, so the
hypothesis~(\ref{eq hyp2}), which is a convenience measure in
the proofs, is not really needed.  Hypothesis~(\ref{eq hyp1}) is
needed to rule out wildly oscillating $f$, since these require different
estimation techniques and Theorem~\ref{th integral test} probably
does not hold for such $f$.  Of course one can prove avoidance for
some such $f$ by comparing to an upper envelope function $\tilde{f} \geq f$.
To see that~(\ref{eq hyp3}) is not too burdensome, note that it
is satisfied in the special cases $f(z) = z^\alpha$ and 
$f(z) = z / \exp (\log^\alpha z)$, which come up naturally in this paper,
and apparently whenever $f''(z)$ behaves in a regular manner.
Note that~(\ref{eq integral test}) is satisfied for $g(z) = 
\exp(\log^{1/2 + \ee} z)$, thus providing a near converse to 
Theorem~\ref{th converse}.  When there
is a gap between first and second moment results, the second moment
result is usually sharp.  Thus Theorem~\ref{th converse} is 
almost certainly not sharp.  But Theorem~\ref{th integral test} is probably
not sharp either, since even
if in principle the second moment method yields a sharp condition via
Lemma~\ref{lem 2mm} below, we do not know whether~(\ref{eq integral test})
is necessary for this.  

\begin{lem} \label{lem reduction}
The probability $p$ of Brownian motion avoiding some set congruent
to $\C_f$ is 0 or 1.  If the probability of Brownian motion avoiding
some $\C_{f , \vv}$ is positive, then $p= 1$. If
for some fixed $\ee \in(0,1)$,
the probability  of Brownian motion avoiding some 
$\C_{(1-\ee)f , \vv}$ is 0, then $p=0$. 
\end{lem}

\noindent{\sc Proof:}  The event that some random set congruent to
$\C_f$ is avoided after a random finite time is a tail event, so
its probability, $p_\infty$, is 0 or 1.  Let $p_t$ be the probability of 
avoiding some random set congruent to $\C_f$ from time $t$ onwards; 
then $p_t \uparrow p_\infty$.  But by the strong Markov property, 
$$p_t = \E_{B(s)} p_{t-s} = p_{t-s} , $$
and therefore $p=p_0 = p_\infty$, proving a zero-one law for
existence of a set congruent to $\C_f$ avoided by Brownian motion.  
Of course it must be 1 if Brownian motion can avoid a $\C_\vv$.
Conversely, if the probability of avoiding some random $\ww + \C_\vv$
is 1, then choosing an arbitrary $\ee\in(0,1)$ and
$\yy \in (\ee / 2) \Z^3$ as close as possible to
$\ww$, there is a positive probability of avoiding some
$\yy + \C_{(1-\ee) f , \vv}$.  By countable additivity this
probability is positive for some fixed $\yy$, and by coupling, for
every fixed $\yy$ and in particular for $\yy = 0$.   $\Cox$
\begin{quote}
{\bf Notation:} Recall that
$$U(L , \theta) = {q(L , \theta) \over q(L)^2} \; .$$
We will identify points in $R^3$ with vectors, in the obvious way.
Let $\theta (\vv , \ww)$ denote the angle between $\vv$ and $\ww$
and let $\theta (\vv) = \theta (\vv , (0,0,1))$.  
\end{quote}

The second moment method is stated in the following lemma.

\begin{lem}[second moment method] \label{lem 2mm}
Suppose that there is a function $U(\theta) \leq \infty$ such that
$$U(L , \theta) \leq U(\theta)$$
for all sufficiently large $L$.  If 
$$\int U(\theta (\vv)) \, dS (\vv) < \infty$$
where $dS$ is surface measure on the unit sphere, then
Brownian motion avoids $f$-thorns.  If furthermore,
$$\int U(\theta (\vv)) \theta^{-\beta} \, dS (\vv) < \infty$$
then the set $ A = \{ \vv : |\vv| = 1 , \cR \cap \C_\vv = \emptyset \}$
of directions 
of $f$-thorns avoided by Brownian motion 
has dimension at least $\beta$, with positive probability.  
\end{lem}

\noindent{\sc Proof:} Let $W_L$ be the measure of the set
$\{ \vv : |\vv| = 1 , \cR \cap \C_\vv \cap \B(0,L) = \emptyset \}$.  
By Fubini's Theorem, $\E W_L = 4 \pi q(L)$.  Another application of 
Fubini's Theorem gives
$$\E W_L^2 = \int \int \P (\tau_L < \tau_{\C_\vv} \wedge \tau_{\C_\ww})
   \, dS(\vv)dS(\ww) = 4 \pi \int q (L , \theta (\vv)) \, dS (\vv) .$$
By Cauchy-Schwartz, $\P (W_L > 0) \geq (\E W_L)^2 / \E W_L^2$,
and this in turn is, up to a factor of $4 \pi$, equal to the reciprocal of 
$\int U(L,\theta (\vv)) \, dS (\vv)$.  
Thus finiteness of 
$\int U(\theta (\vv)) \, dS (\vv)$ implies that $\P (W_L > 0)$
is bounded away from zero for large $L$, which implies that 
$$\P \left ( \lim_{L \rightarrow \infty} \one_{W_L > 0} = 1 
   \right ) > 0 .$$
This and the previous lemma complete the proof of the first statement.  

For the second statement, let $\Xi = \bigcap_{\ee > 0} \Xi_\ee$ 
be a random nonempty subset of the unit 
sphere with the property that if $\xx , \yy$ are points of the 
unit sphere, then
$$\P (\xx \in \Xi_\ee) \geq C \ee^\beta$$ 
and
$$\P (\xx , \yy \in \Xi_\ee) \leq C \ee^{2 \beta} |\xx - \yy|^{- \beta} .$$
It is shown in~\cite[Lemma 5.1]{Peres} how to construct such sets using
a Cantor-like construction, and that any set $A$ with 
$\P (\Xi \cap A \neq \emptyset) > 0$ must have dimension at least $\beta$.  

Construct the sets $\Xi_\ee$ independent of the Brownian motion.  
Recall that
$ A = \{ \vv : |\vv| = 1 , \cR \cap \C_\vv = \emptyset \}$ and let  
$A_L = \{ \vv : |\vv| = 1 , \cR \cap \C_\vv \cap \B(0,L) = \emptyset \}$.
Let $A_L' = A_L \cap \Xi_{1/L}$.  Let $W_L'$ be the measure of $A_L'$.  
Then Fubini's Theorem gives  
$$\E W_L' \geq C q(L) L^{- \beta}$$
and
$$\E (W_L')^2 \leq C \int q(L , \theta (\vv)) L^{- 2 \beta} 
   |\theta (\vv)|^{- \beta} \, dS(\vv) . $$
Thus by Cauchy-Schwartz again, $\P (W_L' > 0)$ is at least a constant times
the reciprocal of $\int U(\theta (\vv )) |\theta (\vv)|^{-\beta} \, dS(\vv)$.
If this integral is finite then 
$$\P (A \cap \Xi \neq \emptyset) = \P (\bigcap_L A_L' \neq 0) \geq
   \lim \sup \P (W_L' > 0) > 0 ,$$
which shows that $A$ intersects $\Xi$ with positive probability, and hence
has dimension at least $\beta$, with positive probability.
$\Cox$

As mentioned before, the estimates on $q(L)$ and $q(L , \theta)$ that
we plug into this lemma in order to prove Theorems~\ref{th 2.3} 
and~\ref{th integral test} are proved in subsequent sections.
We end this section with a proof of Theorem~\ref{th converse}.
The version of the first moment method that we need is:
\begin{lem} \label{lem 1mm}
Fix any $\ee > 0$.  Recall that $g(z) = z / f(z)$ and suppose that 
\begin{equation} \label{eq 1mm}
\lim_{L \rightarrow \infty} q(L) g(L) = 0 .
\end{equation}
Then Brownian motion does not avoid $(1+\ee) f$-thorns.  
\end{lem}

\noindent{\sc Proof:} By Lemma~\ref{lem reduction}, it suffices to show
that Brownian motion avoids no $\C_{(1+\ee/2) f , \vv}$.  On
the event $H$ that Brownian motion avoids some $\C_{(1+\ee/2) f , \vv}$, 
stopping at $\tau_L$ we see that $\tau_L < \tau_{\C_\ww}$ for
every $\ww$ such that $\theta (\ww , \vv) < (\ee/4) f(L) / L = \ee/(4 g(L))$.
Thus on $H$ there is a $c > 0$ such that $W_L \geq c /g(L)$.
Recall that $\E W_L = 4 \pi q(L)$.
This and~(\ref{eq 1mm}) imply $\P (W_L \geq c /g(L)) \rightarrow 0$
for all $c > 0$. Thus, $\P (H) = 0$, which finishes the proof.
$\Cox$

\noindent{\sc Proof of Theorem}~\ref{th converse}: Fix $g(z) = 
\exp (\alpha \log^{1/2} z)$.  Let $L = 2^k$ for some integer $k$.
We compute $q(L)$ as follows.  Define
$$A_j = \{ \vv := (x,y,z) : 2^j < z \mbox{ and } |\vv| \leq 2^{j+1}
   \mbox{ and } \sqrt{x^2 + y^2} \leq f(2^j) \} \, .$$
Then $\C$ contains the disjoint union $A$ of the sets $A_j$ and
the event $\{ \tau_L < \tau_\C \}$ implies the event
$\{ \tau_L < \tau_A\}$.  Conditioning on
successive values of $B(\tau_{2^j})$ and using the strong Markov 
property, we get
$$q(L) \leq \prod_{j=1}^k \P_{B(\tau_{2^j})} \left ( B(s) \notin
   A_j \mbox{ for } 0 < s < \tau_{2^{j+1}} \right ) \, .$$

The terms of the product may be bounded as follows.  Scaling
down by a factor of $2^{j+1}$ transforms $B(\tau_{2^j})$ into a point
on the sphere of radius $1/2$ and $A_j$ into a superset of
a cylinder whose axis is the segment $[1/2 , 3/4]$ and
whose radius is 
$${f(2^j) \over 2^{j+1}} = {1 \over 2} \, \exp \left ( - \alpha 
   \sqrt{j \log 2} \right ) .$$
The probability of avoiding a cylinder of length $1/4$ and
radius $r$ before hitting the boundary of the unit sphere, 
starting at a point of modulus $1/2$, is bounded above by
$1 - K / |\log r|$ for some constant $K$, and since
in our case $|\log r| = (\alpha \sqrt{\log 2} + o(1)) \sqrt{j}$, we get
$$q (L) \leq O(1) \prod_{j=1}^k \left ( 1 - {K_1 \over \alpha \sqrt{j}} 
   \right ) \leq O(1) \exp \left ( - \sum_{j=1}^k 
   {K_1 \over \alpha \sqrt{j}} \right )  \leq \exp (- \sqrt{k}) $$
when $\alpha$ is small.  Thus for sufficiently small $\alpha$,
$q(L) g(L) \leq \exp ((\sqrt{\log 2}\alpha - 1) \sqrt{k}) \rightarrow 0$,
which together with Lemma~\ref{lem 1mm} proves 
that Brownian motion does not avoid $(1+\ee) f$-thorns.
It remains to notice that $(1+\ee) z / \exp (\alpha \log^{1/2} z)
\leq z / \exp ((\alpha - \ee) \log^{1/2} z)$
for large $z$, to complete the proof of
Theorem~\ref{th converse}.

\section{Proof of Theorem~\protect{\ref{th 2.3}}}

\setcounter{equation}{0}

The theorem is proved via the second moment method, following immediately
from the estimate:
\begin{lem} \label{lem 2.3 estimate}
Suppose $f(z) = z^\gamma$.  Then
there are constants $M$ and $\beta$ for which 
\begin{equation} \label{eq 2.3}
U (L , \theta) \leq M |\log \theta|^{1 / (1 - \gamma)} 
  \log^\beta |\log \theta| .
\end{equation}
\end{lem}
Lemma~\ref{lem 2.3 estimate} will be proved at the end
of Subsection 3.3.

Recall that $B(t) = (X(t),Y(t),Z(t))$ and
let $V(t) = \sqrt{X(t)^2 + Y(t)^2}$.  
For any process $\{ \Lambda(t) \}$,
let $T^\Lambda (L) = \inf \{ t > 0 : \Lambda(t) = L \}$ be the time to hit the
value $L$, so that the notation $\tau_L$ is the same as $T^{|B|} (L)$.
For the duration of this section, fix a $\gamma \in (0,1)$
and define $f(z) = z^\gamma$.  Since $f$ is fixed, we suppress it 
from the notation.   Define $\C^{L_0} = \C \setminus (\R^2 \times 
(-L_0 , L_0))$ to be the part of $\C$ with $z$-coordinate at least
$L_0$ in magnitude.  Let $\C_{\theta} = \C_{\vv_\theta}$ as before,
and define $\C^{L_0}_\theta = \{ \yy \in \C_\theta : 
|\yy \cdot \vv_\theta| \geq L_0 \}$ analogously to $\C^{L_0}$.  
Frequent use is made of the following fact (see~\cite[proof
of Theorem 4.3.8, p.~103]{Knight}),
\begin{quote}
{\bf Fact (*):}  For any $0 < v_1 < v_2 < v_3$ and any point
$\vv = (x,y,z)$ such that $x^2 + y^2 = v_2^2$, the hitting
probabilities for the radial process $V(t)$ obey
$$\P_\vv (T^V (v_1) < T^V ({v_3})) = {\log v_3 - \log v_2
   \over \log v_3 - \log v_1} \, .$$
\end{quote}

Fix parameters $\alpha > \beta > 2$, to be used throughout
Section 3 (they are different from $\alpha$ and $\beta$
in other sections).  
The proof of~(\ref{eq 2.3}) is based on estimating 
$q(L)$ and $q(L,\theta)$ separately.  To begin, record the
following useful bound.
\begin{pr} \label{pr 1}
There exists an absolute constant $\cl < 1$ such that if $b \geq 2a > 0$,
$\vv = (x,y,z)$, $|z| \leq a$ and $x^2 + y^2 \leq a^2$, then
$$\P_\vv (T^V (a) \geq T^{|Z|} (b)) 
   \leq \cl^{b/a} .$$
\end{pr}

\noindent{\sc Proof:}  It is elementary to see that there is a $\clo < 1$
independent of $\vv$ such that under the above conditions,
$$\P_\vv (T^V (a) \geq T^Z (z+a) \wedge T^Z (z-a)) \leq \clo .$$
By applying the strong Markov property at the times when 
$B(t)$ hits the planes $\{ (x' , y' , z') : z' = z + j a \}$ for integer
values of $j$, we obtain
$$P_\vv (T^V (a) \geq T^{|Z|} (b)) \leq (\clo)^{b/a - 1}  .$$
Now let $\cl = \sqrt{\clo}$ and observe that $(\clo)^{b/a - 1} \leq
\cl^{b/a}$ as long as $b/a \geq 2$.   $\Cox$

\subsection{Estimating $q(L)$}

Define sequences of constants $m_k = k (\log k)^\alpha$, $r_k =
e^{m_k}$ and $q_k = r_k / (\cll \log k)$, where $\cll$ is chosen 
so that $\cll \log \cl \leq -2$ and $\cl$ is the constant from 
Proposition~\ref{pr 1}.  
\begin{lem} \label{lem P(L)}
Let $j(L)$ be the smallest integer $j$ for which $r_{j-1} \geq L$.
There are constants $k_0$ and $\clll$ for which the following estimate holds.
\begin{equation} \label{eq 7} 
q (L) \geq \clll \exp \left ( - \sum_{k = k_0}^{j(L)} {1 \over 1 - 
   \gamma} \left ( {1 \over k} + {\alpha \over k \log k} \right ) \right ) .
\end{equation}
\end{lem}

\noindent{\sc Proof:}  The constant $k_0$ will be chosen large enough
so that certain inequalities hold; we use the usual convention of
replacing $k_0$ by something larger when necessary to satisfy each
subsequent inequality.  The method of achieving a lower bound
on $q(L)$ is to require something stronger, namely that the radial
part $V(t)$ reach $q_k$ before the $z$-component reaches magnitude
$r_k$ for each $k$.  With hindsight (i.e., comparing to the upper
bound at the end of this section) we can see that this method is
sharp up to a constant factor: conditional on avoiding $\C$ up to time
$\tau_L$ it will be true with probability bounded away form zero that
$V$ reaches each $q_k$ before $|Z|$ reaches $r_k$.  

Suppose that $\vv = (x,y,z)$ with $|z| \leq r_{k-1}$ and 
$x^2 + y^2 = q_{k-1}^2$.  
Then
\begin{equation} \label{eq 1}
\P_\vv (T^V(q_k) < T^{|Z|} (r_k) \wedge T^B (\C)) \geq \P_\vv (T^V (q_k)
   < T^V (r_k^\gamma)) - \P_\vv (T^V (q_k) \geq T^{|Z|} (r_k)) .
\end{equation}
We have
\begin{eqnarray}
&& \P_\vv (T^V (r_k^\gamma) \leq T^V (q_k)) \\[2ex]
& = & {\log q_k - 
   \log q_{k-1} \over \log q_k - \log (r_k^\gamma)} 
   \nonumber \\[2ex]
& = & {k (\log k)^\alpha - \log \log k - (k-1) (\log (k-1))^\alpha + 
   \log \log (k-1) \over (1 - \gamma) k (\log k)^\alpha - \log \cll -
   \log \log k } \; . \label{eq 2}
\end{eqnarray}
Our next goal is to simplify this expression.  First we observe that
\begin{equation} \label {eq 3}
k (\log k)^\alpha - (k-1) (\log (k-1))^\alpha = (\log k)^\alpha
   + (k-1) [(\log k)^\alpha - (\log (k-1))^\alpha] .
\end{equation}
Next we apply the Taylor series expansion.  For $k > k_0$,
\begin{eqnarray} 
&& (\log (k-1))^\alpha \\[2ex]
& \geq & (\log k)^\alpha - {\alpha \over k} 
   (\log k)^{\alpha - 1} + {1 \over (k-1)^2} \left [ (\alpha - 1)
   \alpha (\log (k-1))^{\alpha - 2} - \alpha (\log (k-1))^{\alpha - 1}
   \right ] \nonumber \\[2ex]
& \geq & (\log k)^\alpha - {\alpha \over k} (\log k)^{\alpha - 1} -
   {1 \over k^2} \alpha (\log k)^{\alpha - 1} \label{eq 4} ,
\end{eqnarray} 
since the difference between $k^{-2}$ and $(k-1)^{-2}$ 
is $O(k^{-3}) = o((k-1)^{-2} (\log(k-1))^{\alpha - 2})$.
This combined with~(\ref{eq 3}) gives
\begin{eqnarray*}
k (\log k)^\alpha - (k-1) (\log (k-1))^\alpha & \leq & (\log k)^\alpha
   + (k-1) {\alpha \over k} (\log k)^{\alpha - 1} + (k-1) {1 \over k^2}
   \alpha (\log k)^{\alpha - 1} \\[2ex]
& = & (\log k)^\alpha + \alpha (\log k)^{\alpha - 1} - {1 \over k^2}
   \alpha (\log k)^{\alpha - 1} .
\end{eqnarray*}
Thus for $k > k_0$, throwing out two negative terms 
\begin{equation} \label{eq 5}
\P_\vv (T^V (r_k^\gamma) \leq T^V (q_k)) \leq  {(\log k)^\alpha
   + \alpha (\log k)^{\alpha - 1} + \log \log (k-1) \over (1 - \gamma)
   k (\log k)^\alpha - 2 \log \log k} \; . 
\end{equation}
The following inequality is valid for $a,b,c > 0$ such that $b 
\geq 2c$:
$${a \over b-c} \leq {a \over b} + {2ac \over b^2} .$$
This and~(\ref{eq 5}) imply that for large $k$, 
\begin{eqnarray}
\P_\vv (T^V (r_k^\gamma) \leq T^V (q_k)) & \leq & {(\log k)^\alpha
   + \alpha (\log k)^{\alpha - 1} + \log \log (k-1) \over (1 - \gamma)
   k (\log k)^\alpha} \nonumber \\[1ex]
&& + {4 (\log \log k) ((\log k)^\alpha + \alpha (\log k)^{\alpha - 1} 
   + \log \log (k-1)) \over (1 - \gamma)^2 k^2 (\log k)^{2 \alpha}}
   \nonumber \\[3ex]
& \leq & {1 \over 1 - \gamma} \left ( {1 \over k} + {\alpha \over k \log k}
  + {2 \log \log (k-1) \over k (\log k)^\alpha} \right ) \; . \label{eq 6}
\end{eqnarray}
Recall that $\cll \log \cl \leq -2$ and choose $k_0$ such that for $k 
\geq k_0$ we have $r_{k-1} < q_k$.  For $k \geq k_0$ we can apply
Proposition~\ref{pr 1} to obtain
$$\P_\vv (T^V (q_k) < T^{|Z|} (r_k)) \leq \cl^{r_k / q_k} = 
   \cl^{\cll \log k} \leq k^{-2} .$$
Together with~(\ref{eq 1}) and~(\ref{eq 6}), this yields for large $k$,
$$\P_\vv (T^V (q_k) < T^{|Z|} (r_k) \wedge T^B (\C)) \geq
   1 -  {1 \over 1 - \gamma} \left ( {1 \over k} + {\alpha \over k \log k}
   + {2 \log \log (k-1) \over k (\log k)^\alpha} \right ) - {1 \over k^2} .$$
The $1 / k^2$ term is small enough to be absorbed into the last 
error term, so setting
$$p_k \deq {1 \over 1 - \gamma} \left ( {1 \over k} + {\alpha \over k \log k}
  + {4 \log \log (k-1) \over k (\log k)^\alpha} \right ) ,$$
we have finally,
$$\P_\vv (T^V (q_k) < T^{|Z|} (r_k) \wedge T^B (\C)) \geq
   1 - p_k .$$
It remains to multiply these estimates together qua conditional
probabilities.  

Recall that $j$ is defined so that 
$r_{j-2} < L \leq r_{j-1}$ and that all our estimates are
valid for $k \geq k_0$.  
The strong Markov property applied at each time $T^V (q_k)$ implies that
$$q(L) \geq c'_q\prod_{k = k_0}^j (1 - p_k) .$$
For small $a > 0$ we have $\log (1-a) \geq - a - a^2$ and so
(enlarging $k_0$ if necessary, and thereby introducing a constant factor)
$$\log \left ( \prod_{k=k_0}^j (1 - p_k) \right ) 
   \geq \hat c + \sum_{k=k_0}^j (- p_k - p_k^2) .$$
Since $p_k^2$ is summable, as is the lowest order term
$\log \log k / k (\log k)^\alpha$ in the definition of $p_k$,
we get
$$\log \left ( \prod_{k=k_0}^j (1 - p_k) \right ) 
   \geq \cllo - \sum_{k=k_0}^j {1 \over 1 - \gamma}
   \left ( {1 \over k} + {\alpha \over k \log k} \right ), $$
which finishes the proof of~(\ref{eq 7}).   $\Cox$

\subsection{Estimating $q(L,\theta)$}
 
Now begins the task of estimating the probability $q(L,\theta)$
of avoiding both $\C$ and $\C_{\theta}$ until $\tau_L$.  Let $T_k
= \tau_{r_k}$.  The argument proceeds by estimating the conditional 
probability of avoiding both $\C$ and $\C_{\theta}$ between each 
$T_k$ and $T_{k+1}$ given $B(T_k)$, and multiplying the supremum of
these conditional probabilities to give an upper bound.  For 
values of $k$ greater than some $k_1 (\theta)$, this will
be close to 
$$1 - {2 \over 1 - \gamma} \left ( {1 \over k+1} + {\alpha
   \over (k+1) \log (k+1)} \right ) \; , $$ 
corresponding to intersections with $\C$ and $\C_{\theta}$
being roughly independent, while for small $k$ the $2/(1-\gamma)$
is replaced by a $1/(1-\gamma)$.  Multiplying these together
and identifying the value of $k_1$ will then give an upper
bound on $q(L,\theta)$.  The place this bound loses sharpness is
that $k_1$ must be chosen large enough to give 
a leading term of $2 / ((1-\gamma)k)$ even in the worst case, 
that being the case $B(T_k) \in \C$ which is not likely to happen.  

Again define sequences of constants: $a_k = r_k k (\log k)^\beta$,
$\rho_k = \cll a_k \log k$, $b_k = r_{k+1} / (\cll \log k)$ and
$d_k = r_k \cll \log k$, where $\cll$ still satisfies $\cll \log \cl
\leq -2$.
Set $k_1 (\theta) = 
\exp (|\log \theta|^{1/(\alpha - 1)})$.
Lemma~\ref{lem small k} and Corollary~\ref{cor large k} provide 
estimates for small $k$ and large $k$ respectively. 
\begin{lem} \label{lem small k}
Let 
\begin{equation} \label{eq 15}
t_k \deq  {1 \over 1 - \gamma} \left ( {1 \over k+1} + {\alpha
   \over (k+1) \log (k+1)} \right ) - {\clloo \over k 
   (\log k)^{\beta - 1}} .
\end{equation}
There are constants $k_0$ and $\clloo$ such that for any 
$k \geq k_0$ and any $\vv$ with $|\vv| = r_k$,
$$\P_\vv (B[0,T_{k+1}] \cap \C^{d_k} \neq \emptyset) \geq t_k .$$
\end{lem}
\begin{lem} \label{lem large k}
For $k \geq k_1 (\theta) \deq \lceil \log (1 / \theta) \rceil$ and 
$|\vv| = r_k$, 
\begin{equation} \label{eq 19}
\P_\vv \left ( B[0 , T_{k+1}] \cap \C^{d_k} \neq
   \emptyset \mbox{ and } B[0 , T_{k+1}] \cap \C_{ 
   \theta}^{d_k} \neq \emptyset \right ) \leq {9 \over (1 - \gamma)^2
   k (\log k)^\alpha} \; .
\end{equation}
\end{lem}
\begin{cor} \label{cor large k}
For $k \geq k_1 (\theta)$ and $|\vv| = r_k$, 
\begin{eqnarray} 
&&\P_\vv \left [ B[0,T_{k+1}] \cap (\C^{d_k} \cup \C_{
   \theta}^{d_k}) \neq \emptyset \right ] \nonumber \\
& \geq & {2 \over 1 - \gamma} \left ( {1 \over k+1} + 
  {\alpha \over (k+1) \log (k+1)} - {\clloo \over k (\log k)^{\beta - 1}} 
   \right ) - {9 \over (1 - \gamma)^2 k (\log k)^\alpha}  
   \; . \label{eq s_k}
\end{eqnarray}   $\Cox$
\end{cor}

\noindent{\sc Proof of Lemma}~\ref{lem small k}:  Let
$$S_k = \inf \{ t > 0 : \exists s \in (0,t) : V(s) = a_k 
   \mbox{ and } V(t) = r_k^\gamma \}, \qquad \inf \emptyset = \infty,$$
be the first time when $V(t) = r_k^\gamma$ after the first time
that $a_k$ is hit by $V$.  
Letting $\E_\yy$ denote expectation with respect to $\P_\yy$ we have 
for every $\yy$ of modulus $r_k$,
\begin{eqnarray}
&& \P_\yy (B[0,T_{k+1}] \cap \C^{d_k} \neq \emptyset) 
   \nonumber \\[1ex]
& \geq & \E_\yy \P_{B(T^V(a_k))} (T^V (r_k^\gamma) < T^V (b_k)) -
   \P_\yy (T^V (a_k) \geq T^{|Z|} (\rho_k)) \nonumber \\
&& - \E_\yy \one_{\{|Z(T^V(a_k))| \leq \rho_k\}} \P_{B(T^V(a_k))}
   (T^{|Z|} (r_{k+1}) < T^V (b_k)) - \P_\yy (|Z({S_k})| \leq d_k) . 
   \label{eq A}
\end{eqnarray}
In words, this says: wait until the radial part reaches $a_k$ then
see if it comes back to $r_k^\gamma$ before reaching $b_k$; if it
does, it must hit $\C^{d_k}$ at this time $S_k$  unless the 
$z$-coordinate is wrong.  This is covered by the union of three 
events: (1a) $|Z|$ might reach $\rho_k$ before the radial part 
reaches $a_k$; (1b) $|Z|$ might reach $r_{k+1}$ before $S_k$, despite
having magnitude at most $\rho_k$ at time $T^V (a_k)$; 
or (2) $|Z|$ may be smaller than $d_k$ at time $S_k$.  The point
of waiting for the radial part to reach $a_k$ before coming back
is to make event (2) unlikely.  We give easy estimates on these 
three probabilities before doing the Taylor series computation
for the probability of the radial part coming back to $r_k^\gamma$ before
hitting $b_k$.

For (1a) we use Proposition~\ref{pr 1}.  Recalling that $|\yy| = r_k$
gives
\begin{equation} \label{eq 8}
\P_\yy (T^V(a_k) < T^{|Z|} (\rho_k)) \geq 1 - \cl^{\rho_k / a_k} = 1 - 
   \cl^{\cll \log k} \geq 1 - k^{-2} .
\end{equation}
For (1b), condition on $B(T^V(a_k))$ to get:
\begin{eqnarray*}
&&\E_\yy 1_{\{|Z(T^V(a_k))| \leq \rho_k\}}
   \P_{B(T^V(a_k))} (T^{|Z|} (r_{k+1}) < T^V (b_k)) \\
& \leq & \sup \{ \P_{(x,y,z)} (T^{|Z|} (r_{k+1}) < T^V (b_k)) : 
   x^2 + y^2 = a_k^2 , |z| \leq \rho_k \} .
\end{eqnarray*}
Since $a_k$ and $\rho_k$ are less than $b_k$ for large $k$, 
Proposition~\ref{pr 1} gives
\begin{equation} \label{eq 13}
\E_\yy 1_{\{|Z(T^V(a_k))| \leq \rho_k\}}
   \P_{B(T^V(a_k))} (T^{|Z|} (r_{k+1}) < T^V (b_k)) \leq 
   \cl^{r_{k+1} / b_k} = \cl^{\cll \log k} \leq k^{-2} .
\end{equation}

For (2), let $A_k$ denote the event that $|Z (S_k)| \leq d_k$.
Since $a_k / r_k = k (\log k)^\beta$, the $\P_\yy$ distribution
of $T^V (a_k)$ is stochastically greater than $(r_k k (\log k)^\beta)^2$
times some fixed distribution.  The $\P_\yy$ distribution of $S_k$
is even greater.  Since $|Z|$ is independent of $V$, the $\P_\yy$
density of $|Z (S_k)|$ is bounded by $\clool / (r_k k (\log k)^\beta)$
for some constant $\clool$.  Hence there exist constants $\clolo$ and
$\cloll$ for which
\begin{equation} \label{eq 14}
\P_\yy (A_k) \leq {\clolo d_k \over r_k k (\log k)^\beta} \leq 
   {\cloll \over k (\log k)^{\beta - 1}} \; .
\end{equation}

For the final estimate, we condition on $B(T^V (a_k))$ and use
the fact that the event in question depends only on $V$ to get
$$\E_\yy \P_{B(T^V(a_k))} (T^V (r_k^\gamma) < T^V (b_k)) = 
   {\log b_k - \log a_k \over \log b_k - \log (r_k^\gamma)} \; .$$
Expanding the RHS according to the definitions gives
\begin{eqnarray*}
&&{(k+1) (\log (k+1))^\alpha - \log \cll - \log \log k - k (\log k)^\alpha
   - \log k  - \beta \log \log k \over
   (k+1) (\log (k+1))^\alpha - \log \cll - \log \log k - \gamma k 
   (\log k)^\alpha} \\[3ex]
& \geq & {(k+1) (\log (k+1))^\alpha - k (\log k)^\alpha - 2 \log k \over
   (1 - \gamma) (k+1) (\log (k+1))^\alpha + \gamma ((k+1) 
   (\log (k+1))^\alpha - k (\log k)^\alpha)} \; .
\end{eqnarray*}
For positive $a,b$ and $c$ we always have 
$${a \over b+c} \geq {a \over b} - {ac \over b^2} $$
and so this is at least
\begin{eqnarray}
&& {(k+1) (\log (k+1))^\alpha - k (\log k)^\alpha - 2 \log k \over
   (1 - \gamma) (k+1) (\log (k+1))^\alpha} \nonumber \\[1ex]
&& - \gamma ((k+1) (\log (k+1))^\alpha - k (\log k)^\alpha )
   {(k+1) (\log (k+1))^\alpha - k (\log k)^\alpha - 2 \log k \over
   (1 - \gamma)^2 (k+1)^2 (\log (k+1))^{2 \alpha}} \nonumber \\[3ex]
& \geq & {(k+1) (\log (k+1))^\alpha - k (\log k)^\alpha - 2 \log k \over
   (1 - \gamma) (k+1) (\log (k+1))^\alpha}  \\[1ex]
&&  - {\gamma [(k+1) (\log (k+1))^\alpha - k (\log k)^\alpha ]^2 \over
   (1 - \gamma)^2 (k+1)^2 (\log (k+1))^{2 \alpha}} \; . \label{eq 9} 
\end{eqnarray}
For large $k$, the Taylor series expansion gives 
$$ (\log k)^\alpha \leq (\log (k+1))^\alpha - {\alpha \over k+1} 
   (\log (k+1))^{\alpha - 1} \; .$$
Hence
\begin{eqnarray}
&& (k+1) (\log (k+1))^\alpha - k (\log k)^\alpha \\[1ex]
& \geq & (\log (k+1))^\alpha + k {\alpha \over k+1} 
   (\log (k+1))^{\alpha -1} \nonumber \\[1ex]
& = & (\log (k+1))^\alpha + \alpha (\log (k+1))^{\alpha -1}
   - {\alpha \over k+1} (\log (k+1))^{\alpha -1} . \label{eq 10}
\end{eqnarray}
On the other hand,~(\ref{eq 4}) with $k$ replaced by $k+1$ gives
\begin{eqnarray}
(k+1) (\log (k+1))^\alpha - k (\log k)^\alpha & = & 
   (\log (k+1))^\alpha + k [(\log(k+1))^\alpha - (\log k)^\alpha] 
   \nonumber \\[1ex]
& \leq & (\log (k+1))^\alpha + \alpha ({k \over k+1} + {k \over
   (k+1)^2}) \log (k+1))^{\alpha - 1} \nonumber \\[1ex]
& \leq & 2 (\log (k+1))^\alpha . \label{eq 11}
\end{eqnarray}
Plugging~(\ref{eq 10}) and~(\ref{eq 11}) into~(\ref{eq 9}) gives,
for large $k$,
\begin{eqnarray}
&&\E_\yy \P_{B(T^V(a_k))} (T^V (r_k^\gamma) < T^V (b_k)) \\[2ex]
& \geq & {(\log(k+1))^\alpha + \alpha (\log(k+1))^{\alpha - 1} 
   - {\alpha \over k+1} (\log (k+1))^{\alpha - 1} - 2 \log k \over
   (1 - \gamma) (k+1) (\log (k+1))^\alpha } \nonumber \\[1ex]
&& - {\gamma (2
   (\log(k+1))^\alpha)^2 \over (1 - \gamma)^2 (k+1)^2 (\log (k+1))^{2
   \alpha}} \nonumber \\[2ex]
& \geq & {1 \over 1 - \gamma} \left ( {1 \over k+1} + {\alpha 
   \over (k+1) \log (k+1)} - {4 \log (k+1) \over (k+1) 
   (\log (k+1))^\alpha} \right ) - {4 \gamma \over (1 - \gamma)^2
   (k+1)^2} \nonumber \\[2ex]
& \geq & {1 \over 1 - \gamma} \left ( {1 \over k+1} + {\alpha 
   \over (k+1) \log (k+1)} - {8 \over (k+1) 
   (\log (k+1))^{\alpha - 1}} \right )  . \label{eq 12}
\end{eqnarray}

All the parts of inequality~(\ref{eq A}) have now been estimated.
Plugging in~(\ref{eq 12}),~(\ref{eq 8}),~(\ref{eq 13}) and~(\ref{eq 14})
gives
\begin{eqnarray*}
&& \P_\yy (B[0,T_{k+1}] \cap \C^{d_k} \neq \emptyset) \\[1ex] 
& \geq & {1 \over 1 - \gamma} \left ( {1 \over k+1} + {\alpha
   \over (k+1) \log (k+1)} - {8 \over (k+1) (\log (k+1))^{\alpha - 1}}
   \right ) - {2 \over k^2} - {\cloll \over k (\log k)^{\beta - 1}} 
\end{eqnarray*}
which may be written in the form~(\ref{eq 15}) thus proving
Lemma~\ref{lem small k}.   $\Cox$

\noindent{\sc Proof of Lemma}~\ref{lem large k}:  
If the event in~(\ref{eq 19}) occurs then it occurs as follows: 
$B(t)$ hits one of the two
sets $\C^{d_k}$ or $\C_{ \theta}^{d_k}$, and then the other.  
By symmetry, the probability of this is at most twice 
the supremum of the probability 
of hitting $\C_{\theta}^{d_k}$ and then hitting $\C^{d_k}$,
where the supremum is taken over all starting points
$\yy$ with $|\yy| = r_k$.  
Conditioning on 
$T^V(r_k)$ and on the first point $\zz$ where $B(t)$
hits $\C_{ \theta}$ and using the strong Markov
property shows that the probability in~(\ref{eq 19}) is at most
$2(p_1 +  p_2 p_3)$ where
\begin{eqnarray*}
p_1 & = & \sup \{ \P_\yy (T^V (r_k) > T^{|Z|} (d_k)) : |\yy| = r_k \} ; \\
p_2 & = & \sup \{ \P_{(x,y,z)} (T^V (r_{k+1}^\gamma) < T^V (r_{k+1})) :
   x^2 + y^2 = r_k \} ; \\
p_3 & = & \sup \{ \P_\zz (B[0,T_{k+1}] \cap \C^{d_k}_\theta 
   \ne \emptyset) : \zz \in \C_{\theta} \}
= \sup \{ \P_\zz (B[0,T_{k+1}] \cap \C_\theta 
   \ne \emptyset) : \zz \in \C_{\theta}^{d_k} \}
\\
& \leq &\sup \{ \P_\zz (T^V(r_{k+1}^\gamma) < T^V(r_{k+1})) : \zz 
   \in \C_{\theta}^{d_k} \} .
\end{eqnarray*}

To estimate these three probabilities, we use Proposition~\ref{pr 1}
and the Fact (*) twice.  First, by Proposition~\ref{pr 1}, 
when $|\yy| = r_k$, we have 
\begin{equation} \label{eq 16}
\P_\yy (T^V (r_k) > T^{|Z|} (d_k)) \leq \cl^{d_k / r_k} \leq k^{-2} .
\end{equation}
Secondly, for $x^2 + y^2 = r_k^2$, Fact (*) gives
\begin{eqnarray*}
\P_{(x,y,z)} (T^V (r_{k+1}^\gamma) < T^V (r_{k+1})) & = & {\log r_{k+1}
   - \log r_k \over \log r_{k+1} - \log (r_{k+1}^\gamma)} \\[1ex]
& = & {(k+1) (\log (k+1))^\alpha - k (\log k)^\alpha \over (1 - \gamma)
   (k+1) (\log (k+1))^\alpha} \; .
\end{eqnarray*}
Recalling from~(\ref{eq 11}) that 
$$(k+1) (\log (k+1))^\alpha - k (\log k)^\alpha \leq 2 (\log
   (k+1))^\alpha$$
and assuming $k \geq k_0$ then yields
\begin{equation} \label{eq 17}
\P_{(x,y,z)} (T^V (r_{k+1}^\gamma) < T^V (r_{k+1})) \leq {2 \over (1 
   - \gamma) (k + 1)} \; .
\end{equation}

Lastly, assume $k \geq k_1 (\theta) = \lceil \log (1 / \theta) \rceil$.
Then the distance $\Delta$ between $\C^{d_k}$ and 
$\C_{\theta}^{d_k}$ is at least 
$$\theta r_k - 2 r_k^\gamma \geq r_k (e^{-k} - 2 r_k^{\gamma - 1})
   \geq {r_k \over e^{k+1}} $$
for $k \geq k_0$ where $k_0$ is independent of $\theta$.  If
$k \geq k_1 (\theta)$ and $\zz$ is any point on $\C_{
\theta}^{d_k}$, we have 
\begin{eqnarray}
\P_\zz (T^V (r_{k+1}^\gamma) < T^V (r_{k+1}) ) & \leq & {\log r_{k+1}
   - \log \Delta \over \log r_{k+1} - \log (r_{k+1}^\gamma)} \nonumber \\
& \leq & {(k+1) (\log (k+1))^\alpha - k (\log k)^\alpha +
k + 1 \over (1 - \gamma) (k+1) (\log (k+1))^\alpha} \nonumber \\
& \leq  & {2 (\log (k+1))^\alpha 
 \over (1 - \gamma) (k+1) (\log (k+1))^\alpha}
+ {1 \over (1-\gamma)(\log k)^\alpha} \nonumber\\
&\leq & {2 \over (1-\gamma) (k+1)}
+ {1 \over (1-\gamma) (\log k)^\alpha} \nonumber \\
&\leq& {2 \over  (1-\gamma) (\log k)^\alpha}\; . \label{eq 18}
\end{eqnarray}

Putting together~(\ref{eq 16}),~(\ref{eq 17}) and~(\ref{eq 18}) gives
$$2(p_1 + p_2 p_3) \leq 2 \left ( k^{-2} + {4 \over (1 - \gamma)^2 k (\log
   k)^\alpha} \right )$$ 
which proves Lemma~\ref{lem large k}.   $\Cox$ 

\noindent{\sc Proof of Lemma}~\ref{lem 2.3 estimate}
and {\sc Theorem}~\ref{th 2.3}:  
It remains
to multiply all the conditional probabilities.  Recall the definition 
of $t_k$ as the RHS of~(\ref{eq 15}) and let $s_k = 2 t_k - 9 / 
((1-\gamma)^2 k (\log k)^\alpha)$ be the RHS of~(\ref{eq s_k}).  
When $|\yy| = r_k$ the $\P_\yy$-probability of the event
$$ \{ B[0,T_{k+1}] \cap \C^{d_k} \neq \emptyset \} \cup   
   \{ B[0,T_{k+1}] \cap \C_{\theta}^{d_k} \neq \emptyset \} $$
is bounded below by $t_k$ for any $k \geq k_0$ and by $s_k$ in the
case that $k \geq k_1 (\theta)$.  Let $m$ be such that $r_{m+1} < L
\leq r_{m+2}$, i.e., $m = j-3$ where $j$ is defined in Lemma~\ref{lem
P(L)}.  A repeated application of the strong Markov property at times
$T_k$ then gives 
$$q(L,\theta) \leq \cllol \prod_{k=k_0}^{k_1} (1 - t_k) \prod_{k = 
   k_1 + 1}^m (1 - s_k) \; . $$
For small $a > 0$ we have $\log (1 - a) \leq - a - a^2$ and so
$$\log \left ( \cllol \prod_{k = k_0}^{k_1} (1 - t_k) \times \prod_{k =
   k_1 + 1}^m (1 - s_k) \right ) \leq \clllo + \sum_{k = k_0}^{k_1}
   (- t_k - 2 t_k^2) + \sum_{k = k_1 + 1}^m (- s_k - 2 s_k^2) \; .$$
The series $s_k^2$ and $t_k^2$ are summable, being $O(k^{-2})$, and
the series 
$$\sum_{k=1}^\infty {3 \clloo \over k (\log k)^{\beta - 1}} + {9 \over
   (1 - \gamma)^2 k (\log k)^\alpha}$$
is summable as well, which implies that there is a $\cllll$ for which 
\begin{eqnarray*}
&& \log \left ( \cllol \prod_{k = k_0}^{k_1} (1 - t_k) \times \prod_{k =
   k_1 + 1}^m (1 - s_k) \right ) \\
& \leq & \cllll - \sum_{k = k_0}^{k_1}
   {1 \over 1 - \gamma} \left ( {1 \over k+1} + {\alpha \over (k+1)
   \log (k+1)} \right ) \\
&& - \sum_{k=k_1 + 1}^m {2 \over 1 - \gamma} \left ( {1 \over k+1} + 
   {\alpha \over (k+1) \log (k+1)} \right )  \; .
\end{eqnarray*}
Thus, using the estimate of $q(L)$ in Lemma~\ref{lem P(L)} in the last
step, we have
\begin{eqnarray*}
q(L,\theta) & \leq & \exp \left ( \cllll - \sum_{k = k_0}^{k_1}
   {1 \over 1 - \gamma} \left ( {1 \over k+1} + {\alpha \over (k+1)
   \log (k+1)} \right ) \right. \\
&& \left. - \sum_{k=k_1 + 1}^m {2 \over 1 - \gamma} \left ( {1 \over k+1} + 
   {\alpha \over (k+1) \log (k+1)} \right ) \right ) \\[2ex] 
& = & \exp \left ( \cllll + \sum_{k = k_0}^{k_1}
   {1 \over 1 - \gamma} \left ( {1 \over k+1} + {\alpha \over (k+1)
   \log (k+1)} \right ) \right. \\
&& \left. - \sum_{k=k_0}^m {2 \over 1 - \gamma} \left ( {1 \over k+1} + 
   {\alpha \over (k+1) \log (k+1)} \right ) \right ) \\[2ex] 
& \leq & \cloooo \exp \left ( {1 \over 1 - \gamma} \log k_1 + {\alpha
   \over 1 - \gamma} \log \log k_1 \right. \\
&& \left. - \sum_{k = k_0 - 1}^{j-1} {2 \over 1 - \gamma} 
   \left ( {1 \over k + 1} + {\alpha \over (k+1) \log (k+1)} \right )
   \right ) \\[2ex]
& \leq & \cloolo |\log \theta|^{1 / (1 - \gamma)} (\log |\log
  \theta|)^\cloool q(L)^2  \, .
\end{eqnarray*}
This proves Lemma~\ref{lem 2.3 estimate}.  Since the RHS of~(\ref{eq 2.3})
is integrable over the unit sphere, Theorem~\ref{th 2.3} then follows
from Lemma~\ref{lem 2mm}.   $\Cox$

\section{Avoidance of thorns passing an integral test}
 
\setcounter{equation}{0}

Let $f$ and $g$ be fixed functions satisfying the 
hypotheses~(\ref{eq hyp1})~- (\ref{eq hyp3}) of
of Theorem~\ref{th integral test}, and satisfying the integral
test~(\ref{eq integral test}).  
Recall $U(L,\theta)$ from Section 2.
The remainder of the paper
is devoted to proving 
\begin{equation} \label{eq U2}
U(L , \theta) \leq C_f \theta^{-\xi}
\end{equation}
for some constant $C_f$ and arbitrarily small $\xi > 0$.  This, together 
with the second moment lemma, proves Theorem~\ref{th integral test}.  
In this section we outline the proof of~(\ref{eq U2}).

The idea of the proof is that if $h (\xx) = \P_\xx (\tau_L < \tau_\C)$
solves a Dirichlet problem for $\B_L \setminus \C$ and $h_\theta$
is the analogous function when $\C$ is replaced by $\C_\theta$, then
$h \cdot h_\theta$ ``almost'' solves the Dirichlet problem on
$\B_L \setminus (\C \cup \C_\theta)$; evaluating at the origin, 
$q(L , \theta)$ is almost equal to $q(L)^2$.  The correction
term is the integral of $\grad h \cdot \grad h_\theta$ against
the Green's function for Brownian motion absorbed by $\C \cup \C_\theta$,
as stated with some obfuscation in Lemma~\ref{lem greens} below.
Thus to prove~(\ref{eq U2}), it suffices to get good bounds on
$|\grad h|$ and on the Green's function $G^\theta (0 , \xx)$.
The bounds on $|\grad h|$ are somewhat tedious
to derive, being based on geometric arguments that involve first 
getting bounds on $|h|$, but are reasonably straightforward.  

Bounds on $G$, however, are not straightforward, since if we knew $G$ 
we could solve the problem directly.  One approach
is to use the bound $G^\theta (\xx , \yy)$
by the unrestricted Green's function $|\xx - \yy|^{-1}$.  Not only does
this give reasonable results (under a stronger hypothesis 
than~(\ref{eq integral test})), but it may be bootstrapped to 
give better and better bounds on $G^\theta$.  The (transfinite) 
limit of such bootstrapping is to get an implicit inequality
obeyed by $G^\theta$ and $q(\cdot , \theta)$ in the form of
Lemma~\ref{lem G and q} below.  This together with Lemma~\ref{lem greens}
gives an integral inequality satisfied by $U(\cdot , \theta)$,
Lemma~\ref{lem U bound} below, which leads directly to~(\ref{eq U2}).

That being the conceptual outline, we now state 
a sequence of lemmas, including those mentioned above, which 
form the technical breakdown of the necessary steps.
The first two are merely useful and intuitively obvious principles
which are used repeatedly in the remaining proofs.

We start with a few technical changes to our setup.
First of all, we will give a new meaning
to the symbol $\C^L$, different from that in Section 3.
The change is small and will not confuse a reader who forgets,
so we risk the duplication of notation.
We start with a set $\C \cap \B_L$ and smooth it
in an appriopriate way so that the resulting set
has a $C^2$-boundary. Recall that $\C$ is defined
by a twice differentiable function $f$ but it is truncated
near the origin so that the origin is outside $\C$.
The boundary of each of the two components
of $\C\cap \B_L$ is smooth except for a
circle at each end of this truncated set.
We modify the set $\C\cap \B_L$ to obtain $\C^L$
so that (i) the sets $\C\cap \B_L$ and $\C^L$
may differ only in a neighborhood of radius one
around each of the circles mentioned above; (ii) 
the boundary of $\C^L$ is $C^2$-smooth;
(iii) for large $L<L'$, the sets 
$\C^L \cap \B_{L/2}$ and $\C^{L'} \cap \B_{L/2}$
are identical.

We will also need a new definition similar to that of $q(L)$.
Define $\wqt (L)$ to be the
probability of hitting $\partial \B_L$ before
hitting $\C^{L/2}$ for Brownian motion starting from the
origin. The meaning of $\wqt(L,\theta)$ is derived in an analogous
way: it is the probability of avoiding $\C^{L/2} \cup \C^{L/2}_\theta$
until the hitting time of $\partial \B_L$.  We change the meaning
of $U(L , \theta)$, again so that now 
$$U (L , \theta) \deq {\wqt (L , \theta) \over \wqt (L)^2} ;$$
it is elementary to check that the second moment method
applies equally well when $U$ is defined in terms of $\wqt$ as when
$U$ is defined in terms of $q$, so the substitution is not dangerous
and saves us from a page full of tildes or an unfamiliar letter.

We say that $L$ is a {\em regular value} for $f$ and $\theta$ if
$U(L , \theta) \geq U(L / 4 , \theta)$.  Probably all values are
regular, but in lieu of a proof of that we must consider both
alternatives.

\begin{lem} \label{lem skew}
(i)
Suppose $\mu_r$ is the subprobability
hitting measure on $\partial \B_r$ of Brownian motion on the 
domain $\B_r \setminus \C^\rho$ for some $\rho>0$:
$$\mu_r (A) = \P_0 (B (\tau_{\C^\rho \cup \partial \B_r}) \in A) .$$ 
Then the density 
$${d \mu_L \over dS} (\xx)$$
of $\mu_L$ with respect to area $dS$ on the $L$-sphere is an increasing
function of the angle between $\xx$ and the $z$-axis.   

(ii) Suppose $r_1 < r$ and $\rho>0$. For $\xx$ satisfying $|\xx| = r_1$,
the probability $P_\xx(\tau_r < \tau_{C^\rho})$ is an increasing
function of the angle between $\xx$ and the $z$-axis.
\end{lem}

\begin{cor} \label{cor no varying}
Assume that $g(z) \rightarrow \infty$ as $z \rightarrow \infty$.  Then
$\wqt(2L) / \wqt(L) \rightarrow 1$ as $L \rightarrow \infty$. 
\end{cor}

\begin{lem} \label{lem greens}
Fix $f, L$ and $\theta$.  Let $h_1 (\xx) = \P_\xx (\tau_{\partial \B_L} 
< \tau_{\C^{L/2}})$ be the probability of hitting the $L$-sphere before
$\C^{L/2}$ starting at $\xx$.  Similarly, let $h_2 (\xx) = \P_\xx 
(\tau_{\partial \B_L} < \tau_{\C_{\theta}^{L/2}})$ be the probability
of hitting the $L$-sphere before hitting the rotated cylinder 
$\C_{\theta}^{L/2}$.   
Let $G^\theta$ denote the Green's function for the 
region $\B_L \setminus ( \C^{L/2} \cup \C_{\theta}^{L/2})$.  
Then there is
a constant $r_f$ such that for all regular values of $L \geq r_f$,
$$\wqt(L , \theta) \leq
2 \left [ \wqt(L)^2 + \int_{\B_{L/4} \setminus (\C \cup 
   \C_{\theta})} (\grad h_1 (\xx) \cdot \grad h_2 (\xx)) G^\theta (0 , \xx) 
   \, d\xx \right ] \, .$$
\end{lem}

\begin{lem} \label{lem G and q}
There exist an absolute constant $K$ and a constant $R_f$
depending on $f$, such that for any $\theta$, any 
$L \geq R_f$ and any $\xx$ with $|\xx| \geq R_f$, 
$$G^\theta (0 , \xx) \leq K \wqt(|\xx| / 2 , \theta) |\xx|^{-1} \, .$$
For any values of the parameters, one has the weaker bound
$$G^\theta (0 , \xx) \leq |\xx|^{-1} .$$
\end{lem}

\begin{lem} \label{lem h bounds}
Let $h(u,r) = h_1 (x,y,u)$ for any point $(x,y,u)$ such that $x^2 + y^2 =
r^2$.  Suppose $r \leq z \leq L / 3$.  There is a constant $r_f$
and a $\chbounds > 0$ such that for all $z \geq r_f$,
\begin{equation} \label{eq old 33}
h(z,r) \leq \chbounds {\wqt(L) \over \wqt(z)} {\log (r / f(z)) \over \log g(z)} .
\end{equation}
If $r \geq z$ but the other hypotheses remain the same, then
$$h(z , r) \leq \chbounds {\wqt(L) \over \wqt(r)} .$$
\end{lem}

\begin{lem} \label{lem grad bounds}
Recall that $h(u,r) = h_1 (x,y,u)$ for $(x,y,u)$ such that $x^2 + y^2 =
r^2$.  
Assume that $(x,y,z) \in \B_{L/4}$.
If $r \leq z$ and $ r_f \leq z \leq L/2$ with $r_f$ as in the 
previous lemma, then
\begin{equation} \label{eq grad bound 1}
|\grad h(z,r)| \leq K_f {\wqt(L) \over \wqt(z) } {1 \over r \log g(z)}
\end{equation}
where the constant $K_f$ depends on $f$ but not on $L$.  
If $L/2 \geq r > z $, $r >  r_f$, and
if $\rho$ denotes $\sqrt{z^2 + r^2}$, then we have as well,
\begin{equation} \label{eq grad bound 2}
| \grad h(z,r)| \leq K_f {\wqt(L) \over \rho \wqt(\rho) \log g(\rho)} \, .
\end{equation}
Finally, if $r$ and $z$ are both at most $2r_f$, then 
$|\grad h(z,r)| \leq c \wqt(L) $ where $c$ depends on $f$.
\end{lem}

\begin{lem} \label{lem U bound}
There exist constants $ c_f  > 0 $ and $R_f > 1$ and a function $b(r)$
such that for any $\theta$ and for any $L \geq R/4 \geq R_f$, 
\begin{equation} \label{eq gritty bound}
U (L , \theta) \leq b(R)+ c_f (1 + |\log \theta|) \int_R^L {U(s , \theta) 
   \over s \log^2 g(s)} \, ds .
\end{equation} 
\end{lem}

Lemma~\ref{lem skew} and Corollary~\ref{cor no varying}, 
are proved in the next section; these require little computation.  
Lemma~\ref{lem greens} is proved in the section following.  No
computation is required, but the fact that the 
estimate~(\ref{eq grad bound 1}) for $|\grad h|$ only holds away from 
$\partial \B_L$ forces us to restrict the integral to a smaller ball
and results in some extra estimates.
Lemma~\ref{lem G and q} is also proved in the same section.  
Lemmas~\ref{lem h bounds},~\ref{lem grad bounds} 
and~\ref{lem U bound} are proved in the subsequent, final section.

We conclude this section with a proof of Theorem~\ref{th integral test}
from the above results.  

\noindent{\sc Proof of Theorem}~\ref{th integral test}:  
By assumption, $\int^{\infty^-} (z \log^2 g(z))^{-1} dz < \infty$, so we may
choose $R$ large enough so that $R/4 \geq R_f$ and
$$c_f \int_R^\infty (z \log^2 g(z))^{-1} dz < \xi $$
for $\xi$ arbitrarily small.
By Lemma~\ref{lem U bound}, for any $L \geq R$, $U(L , \theta)$
is bounded above by the value $u_\theta (L)$, where $u_\theta $ solves the 
integral equation
$$u_\theta (x) = b(R) + c_f (1 + |\log \theta|) \int_R^x {u_\theta (s) 
   \over s \log^2 g(s)} \, ds \, .$$
Differentiating, one sees that
$$u_\theta'(x) = c_f (1 + |\log \theta|)
{u_\theta (x)\over x \log^2 g(x)}$$
and hence that 
$$u_\theta(x) =  b(R) \exp \left ( c_f (1 + |\log \theta|) \int_R^x 
   {1 \over s \log^2 g(s)} ds \right ) .$$ 
By the choice of $R$, the integral, $\xi$, may be made arbitrarily small, 
and so 
$$U(L , \theta) \leq u_\theta (L) \leq C_f \theta^{-\xi} ,$$
proving~(\ref{eq U2}).  The function $\theta (\vv)^{-\xi}$ is
integrable over the unit sphere, so the second moment method
finishes the proof that Brownian motion avoides $f$-thorns.  In fact
$\theta (\vv)^{-\xi - \beta}$ is integrable for $\beta$ arbitrarily
close to 2 (by picking $\xi < 2 - \beta$), so the second moment method
shows that the dimension of the set of directions of axes of $f$-thorns
avoided by Brownian motion is greater than $\beta$ for any $\beta < 2$,
proving the dimension result.   $\Cox$

\section{Noncomputational proofs}

\setcounter{equation}{0}

\noindent{\sc Proof of Lemma}~\ref{lem skew}:
We prove only part (i) as (ii) has a similar proof.
We will use a skew-product decomposition.  This is a standard technique, 
so we will limit ourselves to the description of the decomposition.  
See~\cite[Section~7.15]{skew} for more information.  
Let $B_t^*$ be $B_t$ reversed at time $\tau_r$, the time when it
hits the sphere of radius $r$.  In other words, $B_t^* = B(\tau_r - t)$
for $t \in (0 , \tau_r)$.  Let $R_t = |B_t^*|$ denote the modulus
and $A_t = \theta (B_t^*)$ denote the angle with the $z$-axis.  Then
$R_t$ is the time-reversal of a stopped 3-dimensional Bessel process and
$A_t$ is a diffusion $\psi_t$ on the interval $[0,\pi]$ time-changed
according to a clock determined by $R_t$ but otherwise independent of 
$R_t$.  The processes are related by $A_t = \psi_{\beta (t)}$ where 
$\beta (t) = \int_0^t R_u^{-2} du$.  We have $\beta (t) \rightarrow
\infty$ as $t \rightarrow \tau_r$.  

Note that $B_t \in \C^\rho$ for some $t < \tau_r$ if and only if $\psi_s
< D_s$ for some $s < \infty$, where $D_s$ is the maximal angle with the
$z$-axis of any vector in $\C^\rho$ of length $R_{\beta^{-1} (s)}$.  

The hitting distribution on a sphere is uniform for Brownian motion 
starting from its center.  In order to prove the lemma, it will suffice 
to show that the probability of hitting $\C^\rho$ before hitting $\partial
\B_r$ for a Brownian motion starting from the center and conditioned
to exit the sphere at $\xx \in \partial \B_r$ is a decreasing function
of the angle $\theta (\xx)$ that $\xx$ makes with the $z$-axis.  This is
equivalent to proving that the probability of $\{ \psi_s < D_s \}$
is a decreasing function of $\theta (\psi_0)$.  

To show this, we use a coupling argument.  We consider a process
$(\tilde{R} , \psi^1 , \psi^2)$ such that $\tilde{R}$ has the
same distribution as $R$ and such that $\psi^1$ and $\psi^2$
have the same transition probabilities as $\psi$ given $R$.  We
let $\psi_0^1 > \psi_0^2$ and require that if $\psi_t^1 = \psi_t^2$
then $\psi_s^1 = \psi_s^2$ for all $s > t$ (in other words the processes
stay coupled if they meet).  Then clearly $\psi^1_t \geq \psi^2_t$
for all $t$, so the probability that $\psi^1_s < D_s$ for some $s$
is smaller than the probability that $\psi^2_s < D_s$ for some $s$.
$\Cox$

\noindent{\sc Proof of Corollary}~\ref{cor no varying}:
Consider an arbitrary $a <1$. Let $dS$ denote the normalized
surface area measure on $\B(0,L)$
and let $\mu_L$ be defined for $\Lambda\subset \partial \B_L$ by
$\mu_L(\Lambda) = \P_0(B(\tau_L \wedge \tau_{C^{L/2}}) \in \Lambda)$.
Choose a small $\delta>0$ such that
the $S$-measure of 
$A = \{\vv \in \B(0,L) : r(\vv) > \delta L\}$
is greater than $\sqrt{a}$. Then
Lemma~\ref{lem skew} implies that $\mu_L(A) > \sqrt{a} \mu_L(\B(0,L))$.
Note that $\C \cap \B(0,2L)$ is a subset of
$\{ \vv : r(\vv) \leq f(2L)\}$. This and Fact (*) imply that
for $\xx \in A$,
$$\P_\xx (\tau_{2L} < \tau_{\C^L} )
\geq \P_\xx (\tau_{2L} < \tau_\C )
\geq  {\log r(\xx) - \log f(2L) \over \log (2L) - \log f(2L)}
\geq {\log \delta + \log g(2L) \over \log 2 + \log g(2L)}.$$
Since $g(2L)\to\infty$ as $L \to \infty$, we have
$$\P_\xx (\tau_{2L} < \tau_{\C^L} ) \geq \sqrt{a}$$
for $\xx \in A$ and large $L$.
By the strong Markov property and the definition of $\mu_L$,
for large $L$, 
\begin{eqnarray*}
\wqt(2L) & = &\int_{\B(0,L)} \P_\xx (\tau_{2L} < \tau_{\C^L} ) \, d\mu_L (\xx) \\
& \geq & \int_A \P_\xx (\tau_{2L} < \tau_{\C^L} ) \, d\mu_L (\xx) \\
& \geq & \sqrt{a} \mu_L(A) \geq a \mu_L(\B(0,L)) = a \wqt(L).
\end{eqnarray*}
Since $a$ can be chosen arbitrarily close to $1$,
the proof is complete.
$\Cox$

\section{Green's function methods}

\setcounter{equation}{0}

We begin this section with a theorem that is not sufficient for our
purposes, but is a cleaner version of the one we will use.
\begin{th} \label{th clean green}
Let $\C_1$ and $\C_2$ be any two closed regions contained in a ball $\B_L$ of 
radius $L$ centered at the origin $0$.  Assume the origin is in neither.
Suppose $h_i$ is harmonic on $\B_L \setminus \C_i$.
Let $G (\cdot , \cdot)$ denote the Green's 
function for the interior of $\B_L \setminus (\C_1 \cup \C_2)$; 
in other words, if $\tau$ is the exit time from $\B_L \setminus 
(\C_1 \cup \C_2)$, then the expected occupation of a set $A$ 
up to time $\tau$ is 
$$\E_\yy \int_0^\tau \one_A (B(t)) dt = \int_A G (\yy , \xx) \, d\xx \, .$$
Then 
\begin{equation} \label{eq green magic}
\P_0 (\tau_L < \tau_{\C_1 \cup \C_2}) = h_1 (0) h_2 (0) + 
   \int_{\B_L \setminus (\C_1 \cup \C_2)} [\grad h_1 (\xx) \cdot
   \grad h_2 (\xx)] G (0 , \xx) \, d\xx \, ,
\end{equation}
provided the integral is absolutely convergent.
\end{th}

\noindent{\sc Proof:} 
We will write $\lap$ for the Laplacian.
The function $\phi (\yy) = \E_\yy \int_0^\tau f(B(t))
dt$ satisfies $\lap \phi = -2 f$ for any continuous $f$ for which
$\E \int_0^\tau |f (B(t))| dt$ is finite.  
Applying the dominated convergence theorem to the defining equation for $G$
we see that if $f(\xx) = \grad h_1 (\xx) \cdot \grad h_2 (\xx)$ then 
$$\E_\yy \int_0^\tau f (B(t)) \, dt = \int f (\xx) \, 
   G (\yy , \xx )\, d\xx  ,$$
the RHS (and hence the LHS) being absolutely integrable by assumption.
Since $f$ is bounded and continuous, we see that the Laplacian in $\yy$
of $\int f(\xx) G (\yy , \xx) d\xx$ is $-2 \grad h_1 (\yy) \cdot 
   \grad h_2 (\yy)$.
By the product rule for $C^2$ functions, 
$$\lap (h_1 h_2) = h_1 \lap h_2 + h_2 \lap h_1 + 2 \grad h_1 \cdot 
   \grad h_2 ; $$
adding this to the equation 
$$\lap \left ( \int f(\xx) G (\yy , \xx) \, d\xx \right ) = -2 f(\yy)$$
and remarking that $\lap h_1 = \lap h_2 = 0$ shows that 
$\lap \Psi = 0$, where 
$$\Psi (\yy) = h_1 (\yy) h_2 (\yy) + \int_{\B_L \setminus (\C_1 \cup \C_2)}
   [\grad h_1 (\xx) \cdot \grad h_2 (\xx)] G (\yy , \xx) \, d\xx \, .$$

The function $h_1 h_2$ has boundary values $1$ on $\partial \B_L$
and $0$ on $\C_1 \cup \C_2$.  Since $G (\yy , \xx) \rightarrow 0$ as
$\yy \rightarrow \partial (\B_L \setminus (\C_1 \cup \C_2))$,
these are the boundary values of $\Psi$ as well.  This forces
$\Psi (\yy) = \P_\yy (\tau_L < \tau_{\C_1 \cup \C_2})$, 
by the maximum principle (see~\cite[Theorem II.1.8]{Bass}), since
both sides are harmonic with the same boundary conditions.  Setting
$\yy = 0$ proves the theorem.    $\Cox$

We wish to apply this to the case where $\C_1 = \C^{L/2}$ and $\C_2
= \C_{\theta}^{L/2}$, plugging in the bounds on $|\grad h_i|$ from  
Lemma~\ref{lem grad bounds}.  The gradient of $h_i (\xx)$ is difficult 
to control near the boundary of $\B_L$.  
The following lemma allows us to restrict 
attention to $\B_{L/4}$.

\begin{lem} \label{lem substitute}
Let $\mu_{L/4 , \theta}$ be the hitting subprobability measure on 
$\partial \B_{L/4}$ 
defined by
$$\mu_{L/4 , \theta} (A) 
= \P_0 (B (\tau_{\C \cup \C_\theta\cup \partial \B_{L/4}}) \in A) .$$ 
Let $h_1 (\xx)$
and $h_2 (\xx)$ be the probabilities from $\xx$ of hitting 
$\partial \B_L$ before hitting $\C^{L/2}$ and $\C^{L/2}_{\theta}$
respectively.  There is a constant $r_f$ such that for any $\theta$ and
any regular $L \geq r_f$, 
$$\wqt(L , \theta) \leq 2 \int h_1 (\xx) h_2 (\xx) \, 
   d\mu_{L/4 , \theta} (\xx) \, .$$
\end{lem}

\noindent{\sc Proof:} By Corollary~\ref{cor no varying} we may choose 
$r_f$ great enough so that $\wqt(L) \geq .9 \wqt(L/4)$ for all $L \geq r_f$.
When $L$ is regular for $\theta$, it follows that 
$$\wqt(L, \theta) \geq \wqt(L/4 , \theta) {\wqt(L)^2 \over \wqt(L/4)^2}
   \geq .81 \wqt(L/4 , \theta) .$$
An upper bound for $\wqt(L , \theta)$ is the probability of avoiding both
$\C^{L/2}$ and $\C^{L/2}_{\theta}$ until $\tau_{L/4}$ and then avoiding 
$\C^{L/2}$ until time $\tau_L$.  By the Markov property this upper bound
is $\int h_1 (\xx) \, d\mu_{L/4 , \theta} (\xx)$.  
A similar bound holds for $h_2$.  Thus we have:
\begin{eqnarray*}
\wqt(L , \theta) & \leq & \int h_1 (\xx) \, d\mu_{L/4 , \theta} (\xx) \\
\wqt(L , \theta) & \leq & \int h_2 (\xx) \, d\mu_{L/4 , \theta} (\xx) \\
\wqt(L/4 , \theta) & = & \int 1 \, d\mu_{L/4 , \theta} (\xx) .
\end{eqnarray*}
Now use the fact that $x + y - 1 \leq xy$ for $x$ and $y$ in $[0,1]$ to get
that when $L \geq r_f$, 
\begin{eqnarray*}
\wqt(L , \theta) & \leq & 4 \wqt(L , \theta) - 2 \wqt(L/4 , \theta) \\[1ex]
& \leq & 2 \int (h_1 (\xx) + h_2 (\xx) - 1) \, d\mu_{L/4 , \theta} 
   (\xx) \\[1ex]
& \leq & 2 \int h_1 (\xx) h_2 (\xx) \, d\mu_{L/4 , \theta} (\xx) .
\end{eqnarray*}
$\Cox$

\noindent{\sc Proof of Lemma}~\ref{lem greens}: Let $\tau$ be the
hitting time for the set 
$\partial \B_{L/4} \cup \C^{L/2} \cup \C^{L/2}_{\theta}$.
The function $\Psi (\xx) := \E_\xx h_1 (B_\tau) h_2 (B_\tau)$ is harmonic
in the interior of $\B_{L/4} \setminus (\C^{L/2} \cup \C^{L/2}_{\theta})$
with boundary conditions $h_1 h_2$, so by the same argument as in the
proof of Theorem~\ref{th clean green}, 
$$\Psi (\yy) = h_1 (\yy) h_2 (\yy) + \int_{\B_{L/4} \setminus 
(\C^{L/2} \cup \C^{L/2}_{\theta})} 
\left [ \grad h_1 (\xx) \cdot \grad h_2 (\xx) \right ] 
  \wt G^\theta (\yy , \xx) \, d\xx ,$$
where $\wt G^\theta (\cdot , \cdot)$ is the Green's function for
$\B_{L/4} \setminus (\C^{L/2} \cup \C^{L/2}_\theta)$.
Since $\wt G^\theta (\cdot , \cdot)$ is less than
$ G^\theta (\cdot , \cdot)$, the Green's function for
$\B_L \setminus (\C^{L/2} \cup \C^{L/2}_\theta)$, we obtain an upper bound
for $\Psi (\yy)$ by replacing $\wt G^\theta$
with $ G^\theta$ in the last formula.
It can be shown just like in the last part of
Lemma~\ref{lem grad bounds} that the gradients are bounded on 
$\B_{L/4} \setminus (\C^{L/2} \cup \C^{L/2}_{\theta})$ 
so that there is no problem with convergence of the integral
(we do not make any assertion about the size of the bound
at this point).
By the previous lemma, $\wqt(L , \theta) \leq 2 \Psi (0)$, which, together
with the formula for $\Psi$ and the fact $h_1 (0) = h_2 (0) = \wqt(L)$
proves the lemma.   $\Cox$

\noindent{\sc Proof of Lemma}~\ref{lem G and q}: The weaker bound comes
from bounding $G^\theta$ by the Green's function $G(\cdot,\cdot)$
for all of $\R^3$.  

Let $S$ be the normalized surface area measure on $\partial \B_L$.
To prove the stronger bound, we first claim that for any fixed $\theta$
and some absolute constant $K$,
$${d\mu_{L , \theta} \over dS} \leq K \wqt(L / 2 , \theta),$$
where $\mu_{L , \theta}$ is the hitting (subprobability) measure 
on $\partial \B_L$ of Brownian motion started at 0 and killed at 
$\partial \B_L \cup \C^{L/2} \cup \C^{L/2}_\theta$.  
Let $\nu_\xx$ denote the $\P_\xx$ law of $B(\tau_L)$, that is, the hitting
distribution on $\partial \B_L$ of an unkilled Brownian motion 
started at $\xx$.
The Harnack principle shows that the densities 
$d\nu_\xx / dS$ are bounded for $\xx \in \B_{L/2}$ by an
absolute constant $K$. Thus 
$$\mu_{L , \theta} (A) \leq \int \nu_\xx (A) \, d\mu_{L/2 , \theta} (\xx) 
\leq K ||\mu_{L/2 , \theta}|| \int_A dS = K \wqt(L/2 , \theta) S(A),$$
proving the claim.  

Now let $A$ be any set disjoint from the ball $\B_r \subseteq \B_L$.  
Letting $\tau$ be the hitting time of 
$\partial \B_L \cup \C^{L/2} \cup 
\C^{L/2}_{\theta}$ and using the strong Markov property at time 
$\tau_r$ gives
\begin{eqnarray*}
\int_A G^\theta (0 , \yy) d\yy & = & \int \left ( \int \int \one_{B_t \in A} 
   \one_{\tau > t} dt d\P_\yy \right ) \, d\mu_{r , \theta} (\yy) \\[1ex]
& \leq & \int \left ( \int \int \one_{B_t \in A} 
   dt d\P_\yy \right ) \, d\mu_{r , \theta} (\yy) \\[1ex]
& \leq & K \wqt(r/2 , \theta) \int \left ( \int \int \one_{B_t \in A} 
   dt d\P_\yy \right ) \, dS_r (\yy) ,
\end{eqnarray*}
by the above claim for $L=r$, where $S_r$ is normalized surface measure 
on $\partial \B_r$.  But this last quantity is just 
$$K \int_A \wqt(r/2 , \theta)  G(0,\yy) d\yy 
\leq K \int_A \wqt(r/2 , \theta)  |\yy|^{-1} d\yy .$$
Letting $A$ shrink around $\xx$ and leting $r \uparrow |\xx|$
then proves the lemma.
$\Cox$

\section{Geometric bounds}

\setcounter{equation}{0}

The following lemma is needed in the proof of 
Lemma~\ref{lem h bounds}. It is a version of the
boundary Harnack principle but we could not find
a version of that theorem that would apply directly
in our case.

\begin{lem} \label{lem kb1}
Suppose that for some $z_0$,
$$ A_1 = \{ (x,y,z) : x^2 + y^2 < c_1^2 z_0^2  , 
c_2 z_0  < z < c_3 z_0  \},$$
$$ A_2 = \{ (x,y,z) : x^2 + y^2 < c_4^2 z_0^2  , 
c_2 z_0  < z < c_3 z_0  \},$$
$$ W = \{ (x,y,z) : x^2 + y^2 = c_1^2 z_0^2  , 
c_2 z_0  < z < c_3 z_0  \},$$
and $\vv = (x_1,y_1,z_1)$ is a point with $c_5 < z_1 < c_6$,
$c_4^2 z_0^2 < x_1^2 + y_1 ^2 <  c_1^2 z_0^2$.
Assume that $ c_4 < c_1$ and $c_2 < c_5 < c_6 < c_3$.
Then there exists $c_7>0$ which depends on 
$c_1, c_2, c_3, c_5 $ and $c_6$ but does not depend on
$c_4 $ or $z_0$, and such that
$$\P_\vv (B(\tau_{A_1^c}) \in W \mid \tau_{A_1^c} < \tau_{A_2})
> c_7.$$
\end{lem}

\noindent{\sc Proof}:
We will prove the lemma for $z_0=1$. The general case follows
by scaling. We will also assume that $c_4 < c_1/4$. The other case
requires minor modifications.

Let $c_8 = \max(c_3 - c_6, c_5 - c_2)$ and
choose $c_9$ so that $\sum_{k=1}^\infty c_9 k 2^{-k} < c_8/2$.
Let $m_k = \sum_{j=k}^\infty c_9 j 2^{-j}$,
$$D_k = \{(x,y,z): \sqrt{ x^2 + y^2} < c_1  2^{-k},
c_5 - m_k < z < c_6 + m_k\},$$
$$W_k = \{(x,y,z) \in \partial D_k: \sqrt{ x^2 + y^2} = c_1  2^{-k}
\},$$
$$U_k = \partial D_k \setminus W_k.$$
Note that for any $\ww \in W_k$, the distance from $\ww$
to $W_{k-1}$ is $c_1 2^{-k}$ but the distance to
$U_{k-1}$ is not less than $c_9 k 2^{-k}$.
It easily follows from Proposition~\ref{pr 1}
that if $\ww \in W_k$ then
$$\P_\ww( \tau_{U_{k-1}} < \tau_{W_{k-1}}) < e^{-c_{10} k}.$$
If $c_1 2^{-k} \geq 2 c_4$ then it is easy to see that
for any $\ww \in W_k$,
$$\P_\ww (\tau _{A_2} < \tau_{\partial D_k}) < c_{11}<1.$$
Hence, for $\ww \in W_k$, assuming $c_1 2^{-k} \geq 2 c_4$,
$$\P_\ww( \tau_{U_{k-1}} < \tau_{W_{k-1}} \mid 
\tau_{\partial D_k} < \tau _{A_2}) < e^{-c_{12} k}.$$
Now suppose that $\vv \in W_n$ where $n$ is the smallest
number such that $c_1 2^{-n} \geq 2 c_4$. If we condition
Brownian motion not to hit $A_2$ between the first hitting
times of $\partial D_k$ and $\partial D_{k-1}$
for $k\leq n$, then, using the strong
Markov property, we see that for such conditioned process
we may have $B(\tau_{\partial D_k} ) \in W_k$ for all
$k = n, n-1, \dots, 2$, with probability not less than
$$\prod_{k=2}^n (1 - e^{-c_{12} k})
\geq \prod_{k=2}^\infty (1 - e^{-c_{12} k}) = c_{13} >0.$$
Hence, Brownian motion conditioned to avoid $A_2$
before exiting $D_1$ can hit $W_1$ with probability
greater than $c_{13}$.
Brownian motion starting from a point of $W_1$
can hit $W$ before hitting any other part of the boundary
of $A_1$ or $A_2$ with probability greater than $c_{14}$,
independent of $c_4$.
An application of the strong Markov property at the hitting time
of $W_1$ shows that 
$$\P_\vv (B(\tau_{A_1^c}) \in W \mid \tau_{A_1^c} < \tau_{A_2})
> c_{13}c_{14}.$$

The same proof applies to $\vv \in W_k$ for $k< n$.
The result can be extended to all points $\vv = (x_1,y_1,z_1)$ 
with $x_1^2+y_1^2 > c_1 2^{-n}$ and
$c_5 < z_1 < c_6$ using the Harnack inequality.
Finally, it extends to $\vv $ with
$c_4 < \sqrt{x_1^2+y_1^2} < c_1 2^{-n}$ by the boundary Harnack
principle. See~\cite{Bass} for the Harnack inequality and
the boundary Harnack principle.
$\Cox$

\noindent{\sc Proof of Lemma}~\ref{lem h bounds}:
The idea is that escaping from $(z , r)$ to $\partial \B_L$ takes
two steps.  First, one has to escape to $\{ r \approx z \}$.
Approximating $\C$ by a cylinder of radius $f(z)$ about the $z$-axis, we 
see that this probability is roughly $\log (r / f(z)) / \log (z / f (z))$.
Secondly, one must escape from radius roughly $z$ to radius $L$.
This probability is the conditional probability of escaping
to radius $L$ given having escaped to radius $z$, and is thus
roughly $\wqt(L) / \wqt(z)$.  When $r < 2 f(z)$, the cylinder approximation
is too course and we need a third step, namely first escaping to
$\{ r \approx 2 f(z) \}$.  We now rigorize this.

Let $\vv_0$ be the point $(x_0, y_0, z_0)$, where all coordinates
are assumed w.l.o.g.\ to be positive.  Let
$r_0 = \sqrt{x_0^2 + y_0^2}$.  We use the unsubscripted symbols
$x,y,z,r$ to refer to the $x,y,z$ and $\sqrt{x^2 + y^2}$ coordinate
functions respectively.  Assume first that $r_0 < 2 f(z_0)$.  
Let $A_0$ be the circle in the $x,z$-plane, centered on the 
$z$-axis and tangent to the curve $x = f(z)$ at the point 
$(z_0 , f(z_0))$.  By assumption~(\ref{eq hyp3}), 
$A_0$ lies completely inside the region $x \leq f(z)$.  Rotating
this circle around the $z$-axis gives a sphere, $A_1$, tangent to
$\partial \C$ at all points with $r = f(z_0)$ and $z = z_0$, 
and lying inside $\C$.  Let $A_2$ be
the sphere with 4 times the radius and the same center.  Clearly
we may write
\begin{eqnarray} 
h_1 (\vv_0) & \leq & \P_{\vv_0} (\tau_{A_2} < \tau_{A_1}) \label{eq part1} \\
&& \times \sup_{\vv \in A_2} \P_\vv (\tau_{(5 \sqrt{2}/4) z_0} < 
   \tau_{\C}) \label{eq part2} \\
&& \times \sup_{\vv \in \partial \B_{(5 \sqrt{2}/4) z_0}} \P_\vv
   (\tau_L < \tau_{\C}) \label{eq part3} .
\end{eqnarray}

To estimate the term~(\ref{eq part1}), use the facts that the radius
$R_1$ of $A_1$ is at least $f(z_0)$ and that the distance $d (\vv_0 , A_1)$ 
from $\vv_0$ to $A_1$ is at most $r_0 - f(z_0)$, to get
$$\P_{\vv_0} (\tau_{A_2} < \tau_{A_1}) = {4 \over 3} \left ( 1 - 
   {R_1 \over d(\vv_0 , A_1)} \right ) \leq {4 \over 3} \left ( 1 - 
   {R_1 \over R_1 + r_0 - f(z_0)} \right ) \leq {4 \over 3}
   {r_0 - f(z_0) \over f(z_0)} $$
$$ \leq {4 \over 3} c_1 \log \left ( 1 + {r_0 - f(z_0) \over f(z_0)}
\right ) = {4 \over 3} c_1 \log (r_0 / f(z_0)).$$

Let $A_3$ be the cylinder $\{ r \leq f (z_0) / 2 \}$.  Let $A_4$ be the
cylinder with radius $5 z_0 / 4$ whose axis is the subinterval
$[3 z_0 / 4 , 5 z_0 / 4]$ of the $z$-axis.  Observe that
$A_4$ lies inside $\B_{5 \sqrt{2} z_0 / 4}$, and that $A_3 \cap
A_4$ lies inside $\C \cap A_4$ (since $f(z_0) / 2 \leq f(z_0/2)
\leq f(3 z_0 / 4)$).  Thus~(\ref{eq part2}) may be bounded above by
$\P_\vv (\tau_{A_4} < \tau_{A_3})$.  When $z_0$ is sufficiently large,
$4 f(z_0) \leq z_0 / 10$, and thus the $z$-coordinate of $\vv$
is in $[.9 z_0 , 1.1 z_0]$ for every $\vv \in A_4$.  By scaling, there
is a uniform lower bound $\ee_1 > 0$ for the probability 
$$\P_\vv (B(\tau_{A_4}) \in \R^2 \times (3z_0 / 4 , 5 z_0 / 4))$$
that a Brownian motion started at $\vv$ exits $A_4$ along the curved
boundary $W := \{ (x,y,z) : x^2 + y^2 = 25 z_0^2 / 16 , 
3z_0 / 4 < z < 5z_0 / 4 \}$.  
By Lemma~\ref{lem kb1}, it follows that
$$\P_\vv (B(\tau_{A_4}) \in W \| \tau_{A_4} < \tau_{A_3}) \geq \ee .$$
Thus the term~(\ref{eq part2}) is at most
$$\ee^{-1} \sup_{\vv \in A_2} \P_\vv (\tau_W < \tau_{A_3}) .$$
Using Fact (*), this gives an upper bound of
$$ \ee^{-1} {\log (4 f(z_0)) - \log (f(z_0) / 2) \over \log (5 z_0 / 4) 
   - \log (f(z_0) / 2)} = \ee^{-1} {\log 8 \over \log (5 g(z_0) / 2)}.$$

To estimate~(\ref{eq part3}), note that by 
Lemma~\ref{lem skew} (ii) the
supremum is achieved at points $(x,y,0)$ such that $x^2 + y^2 = 
25 z_0^2 / 8$.  Let $\vv$ be such a point.  The sphere 
of radius $z_0 / 2$ around $\vv$ is disjoint from $\C$, so applying
the Harnack principle to points $\ww$ in the set $\Lambda$ 
of points on $\partial \B_{5 z_0 / 4}$ 
within distance $z_0 / 4$ from the $x,y$-plane, 
we see that there is a universal constant $C$ such that 
$$\P_\vv (\tau_L < \tau_{\C^{L/2}}) 
\leq C \P_\ww (\tau_L < \tau_{\C^{L/2}}) .$$
Lemma~\ref{lem skew} (i) implies that the $\mu_L$-measure 
of the set of points on $\B_L$ which form an angle
greater than $\psi$ with the $z$-axis is greater than
the normalized surface measure of the same set. Hence,
$$\P_0 (B(\tau_{5 \sqrt{2}z_0 / 4}) \in \Lambda \| 
\tau_{5 \sqrt{2}z_0 / 4} \leq \tau_{\C^{L/2}}) 
   \geq {|\Lambda| \over 4 \pi (5\sqrt{2}z_0/4)^2} = \wt c.$$
Thus
\begin{eqnarray*}
{\wqt(L) \over \wqt(5 \sqrt{2}z_0 / 4)} & = & \P_0 (\tau_L < \tau_{\C^{L/2}} \| 
   \tau_{5 \sqrt{2}z_0 / 4} < \tau_{\C^{L/2}}) \\[1ex]
& \geq & \wt c \, \P_0 (\tau_L < \tau_{\C^{L/2}} \| 
   \tau_\Lambda < \tau_{\C^{L/2}}) \\[1ex]
& = & \wt c \, \E [\P_{\tau_\Lambda} (\tau_L < \tau_{\C^{L/2}}) 
\mid \tau_\Lambda < \tau_{\C^{L/2}}]\\[1ex]
& \geq & {\wt c \over C} \P_\vv (\tau_L < \tau_{\C^{L/2}}) .
\end{eqnarray*}
Thus~(\ref{eq part3}) is bounded above by
$$ (C/\wt c) {\wqt(L) \over \wqt(5 \sqrt{2}z_0 / 4)} 
\leq 2( C/\wt c) {\wqt(L) \over \wqt(z_0)}$$
for $z_0$ sufficiently large.  

Combine the three pieces~(\ref{eq part1})~-~(\ref{eq part3})
to yield the bound in the lemma.

In the case where $z_0 \geq r_0 \geq 2 f(z_0)$, skip the first
step, writing $h_1 (\vv_0)$ as at most 
$$\ee^{-1} \P_{\vv_0} (\tau_W < \tau_{A_3})$$
times~(\ref{eq part3}).  
Fact (*) then gives an upper bound of 
$$h_1 (\vv_0) \leq \ee^{-1} 
{\log  r_0 - \log ( f(z_0)/2) \over \log (5 z_0/4) - \log( f(z_0) / 2)}
  \sup_{\vv \in \partial \B_{(5 \sqrt{2}/4) z_0}} \P_\vv
   (\tau_L < \tau_{\C^{L/2}}) $$
which is at most a constant multiple of 
$${\log (r_0 / f(z_0)) \over \log g(z_0)} 
  \sup_{\vv \in \partial \B_{(5 \sqrt{2}/4) z_0}} \P_\vv
   (\tau_L < \tau_{C^{L/2}}) 
\leq \widehat c\, {\log (r_0 / f(z_0)) \over \log g(z_0)}
 {\wqt(L) \over \wqt(z_0)}
,$$
since $r_0 / f(z_0) \geq 2$.    

Finally, in the case where $r_0 \geq z_0$, we have from
Lemma~\ref{lem skew} (ii) that $h(z_0 , r_0) \leq h(\lambda , \lambda)$,
where $\lambda = \sqrt{(z_0^2 + r_0^2) / 2}$.  Now apply~(\ref{eq old 33})
and observe that by Corollary~\ref{cor no varying}
$\wqt(r_0) = (1 +o(1)) \wqt(\lambda)$ for large $r_f$.
$\Cox$

\noindent{\sc Proof of Lemma}~\ref{lem grad bounds}:
We start with a proof of~(\ref{eq grad bound 1}).  It will suffice to
show that
$$|h(\vv_0) - h(\vv_1)| \leq K_f \dd {\wqt(L) \over \wqt(z_0)} 
   {1 \over r_0 \log g(z_0)}$$
where $\vv_i = (x_i , y_i , z_i)$, $r_i = \sqrt{x_i^2 + y_i^2}$
and $\dd := |\vv_0 - \vv_1|$ is small.  We let $\K$ be the plane such
that $\vv_0$ and $\vv_1$ are symmetric with respect to $\K$.

Consider first the case when $r_0 > 2 f(z_0)$.  Assume 
$\dd < r_0 / 100$.  For $k \geq 1$, define the following regions:
\begin{quote}
Let $S_k = \partial \B (\vv_0 , 2^k r_0 / 8)$. \\[1ex]
Let $A_0$ be the closure of $ \B (\vv_0 , 2 r_0 / 8)$. \\[1ex]
Let $A_k$ be the closure of the spherical shell between 
$S_k$ and $S_{k+1}$ for $k\geq 1$.  \\[1ex]
Let $\wt \C$ be the set symmetric to $\C$ with respect to $\K$.\\[1ex]
Let $D_k = (\C \cup \wt \C) \cap A_k$.  \\[1ex]
\end{quote}
Suppose $\wt \tau$ is a stopping time such that
$\wt \tau \leq \tau_{\C^{L/2}} \wedge \tau_L$ a.s.
Since $h$ is harmonic on $\B_L \setminus \C^{L/2}$, we 
have $h(\xx) = \E_\xx h(B(\wt \tau))$
for $\xx \in \B (\vv_0 , r_0/4)$. 
Couple Brownian motions $B_0$ and $B_1$ started from points 
$\vv_0$ and $\vv_1$ so that they are mirror images in $\K$.  
We will use superscripts 
to denote hitting times for $B_i$.
Let 
$$\wt \tau = \tau^0_\K \wedge \tau^0_L \wedge \tau^0_{\C^{L/2}}
\wedge \tau^1_\K \wedge \tau^1_L \wedge \tau^1_{\C^{L/2}}.$$
Note that $h(B_0(\wt \tau)) = h(B_0(\wt \tau))$ 
if $\wt \tau = \tau^0_\K \wedge \tau^1_\K$. Using
$\E$ for the law of the coupling we then have
\begin{eqnarray}
h(\vv_0) - h(\vv_1) & = & \E \left [ h ( B_0 (\wt \tau) )
   - h (B_1 (\wt \tau)) \right ] \nonumber \\[1ex]
& = & \E \left [ \one_{\tau^0_L \wedge \tau^1_L = \wt \tau}
[h ( B_0 (\wt \tau) )
   - h (B_1 (\wt \tau))] \right ] \label{eq hit L} \\[1ex]
& + & \E \left [ h (B_0 (\wt \tau))  
   \one_{\tau^1_{\C^{L/2}} = \wt \tau} - 
   h (B_1 (\wt \tau))  
   \one_{\tau^0_{\C^{L/2}} = \wt\tau} \right ] \label{eq hit C}.
\end{eqnarray}
We take care of the term~(\ref{eq hit L}) first.
We may bound it above by
$$ \E \one_{\tau^0_L = \wt \tau} 
(1 - h ( B_{1} ( \wt \tau))) =
\int (1 - \P_\xx (\tau_L < \tau_{\C^{L/2}})) \, d \pi (\xx)$$
where $\pi$ is the subprobability measure corresponding to the location
of $B_{1} (\tau^0_L)$ restricted to the event 
$\{ \tau^0_L = \wt\tau\}$.  
This event is contained in
$\{ \tau^0_L < \tau^0_\K\}$ and so
it is clear that the total mass $||\pi||$ of $\pi$ 
is at most a constant multiple of $\dd / L$, since $r_0 , z_0 \leq L/4$.
Comparing $\C^{L/2}$ to the infinite cylinder of radius $f(L)$, 
and the ball $\B_L$ to the analogous cylinder of radius $L$
one sees from Fact (*) that for $\xx \in \B_{5L/8}$ with $z(\xx) = 0$,
$$\P_\xx (\tau_L > \tau_{\C^{L/2}})
\leq { \log L - \log (L/2) \over \log L - \log f(L)}
= {\log 2 \over \log g(L)}.$$
By the Harnack principle applied in the shell 
between $\partial \B_{9L/16}$ and $\partial \B_L$,
$$\P_\xx (\tau_L > \tau_{\C^{L/2}})
\leq {c_2 \over \log g(L)},$$
for all $\xx \in \partial B_{5L/8}$. By the maximum principle,
the same inequality holds for all $\xx$ with $|\xx| \geq 5L/8$.
Since $\vv_i \in \B_{5L/16}$, we have $|B_{1}(\tau^0_L)| \geq 5L/8$
and thus the term~(\ref{eq hit L}) is at most
$$\int [1 - \P_\xx (\tau_L < \tau_{\C^{L/2}})] \, d \pi  (\xx)
\leq  {c_2 \over \log g(L)} ||\pi||
\leq  {c_3 \delta\over L\log g(L)}.$$
Recalling from Lemma~\ref{lem skew} that $\wqt(2x) / \wqt(x) \rightarrow 1$,
it follows easily that $L \wqt(L) \geq c z_0 \wqt(z_0)$ and hence that
$$c_3 {\dd \over L \log g(L)} \leq c' \dd {\wqt(L) \over \wqt(z_0) z_0 \log
   (g(z_0))} \, .$$
Since $z_0 \geq r_0$, this shows that the term~(\ref{eq hit L}) is
bounded by an expression of the form
$$c \dd {\wqt(L) \over \wqt(z_0) r_0 \log (g(z_0))} \, .$$

We now turn to the term~(\ref{eq hit C}).  
The event $\{\tau^1_{\C^{L/2}} = \wt \tau\}$ 
is contained in the union of events
$\{ \tau^1_{D_k} \leq \tau^1_\K \}$
for $k \geq 0$.  
Hence, (\ref{eq hit C}) is bounded
above by
\begin{equation} \label{eq k breakdown}
\E  \left [ h (B_0 (\wt \tau))  
   \one_{\tau^1_{\C^{L/2}} = \wt \tau}\right ]
= \sum_{k \geq 0} \P_{\vv_0} 
(\tau_{D_k} \leq \tau_\K )
   \sup_{x \in D_k} h(x) 
\leq \sum_{k \geq 0} \P_{\vv_0} 
(\tau_{D_k} \leq \tau_\K )
   \sup_{x \in A_k} h(x)
\, .
\end{equation}
We will need the following lemma. Its proof is given at the
end of this section.

\begin{lem} \label{lem final pieces}
Let $p_1$ be the probability that Brownian motion started from 
$\vv_0$ will hit $\B (\vv_0 , r_0/4)$ before hitting $\K$.  Let $p_2$
be an upper bound for the probability that a Brownian motion
starting from a point $\yy \in S_k$ will hit $S_{k+1}$ before $\K$.
Let $p_3$ be the probability that a Brownian motion starting
from a point $\yy \in S_k$ will hit $D_{k+1}$ before $\K$.  Then
$p_2 < 1$ and there are constants $c_i > 0$ and $\alpha < p_2^{-1}$
depending only on $f$ and such that for 
$2f(z_0) \leq r_0 \leq z_0 \leq L$
and $z_0 \geq r_f$, 
\begin{eqnarray}
\sup_{\xx \in A_k} h(\xx) & \leq & c_{4} \alpha^k {\wqt(L) \over \wqt(z_0)}
   {\log (r_0 / f(z_0)) \over \log (g(z_0))} \label{eq n1} \\[1ex]
p_1 & \leq & {c_{5} \dd \over r_0} \label{eq n2} \\[1ex]
p_3 & \leq & {c_{6} \over \log (r_0 / f(z_0))} \label{eq n3}.
\end{eqnarray}
In the case when $r_f \leq z_0 \leq r_0 \leq L$
the estimates~(\ref{eq n1}) and (\ref{eq n3}) are replaced by
\begin{eqnarray}
\sup_{\xx \in A_k} h(\xx) & \leq & c_{7} \alpha^k 
{r_0\over \rho} {\wqt(L) \over \wqt(\rho)}
   {\log g (r_0) \over \log (g(\rho))}\label{eq n4} \\[1ex]
p_3 & \leq & {c_{8} \over \log g(r_0)} \label{eq n5}.
\end{eqnarray}
where $\rho = \sqrt{r_0^2 +z_0^2}$.
\end{lem}

Note that (\ref{eq n3}) provides also an upper
bound for the probablity of hitting $D_0 \cup D_1$
before hitting $\K$ for Brownian motion starting from
$\vv_0$ (see the proof of Lemma~\ref{lem final pieces}
at the end of this section).
Assuming this lemma for the moment, use the strong Markov property
to see that for $k\geq 0$ and $2f(z_0) \leq r_0 \leq z_0$, 
$$\P_{\vv_0} ( \tau_{D_k} \leq \tau_\K)
\leq p_1 p_3 p_2^k .$$
Combining this with (\ref{eq k breakdown}) and
(\ref{eq n1}) gives
$$
\E  \left [ h (B_0 (\wt \tau))  
   \one_{\tau^1_{\C^{L/2}} = \wt \tau}\right ]
\leq \sum_{k \geq 0} c p_1 p_2^k p_3 \alpha^k {\wqt(L) \over \wqt(z_0)}
   {\log (r_0 / f(z_0)) \over \log (g(z_0))}$$
which reduces to
$$c {\dd \over r_0} {\wqt(L) \over \wqt(z_0) \log (g(z_0))} ,$$
and finishes the proof in the case $2f(z_0) \leq r_0 \leq z_0$.  

The proof of~(\ref{eq grad bound 2}) is completely analogous ---
estimates~(\ref{eq n4}) and~(\ref{eq n5}) have
to be used in place of~(\ref{eq n1}) and~(\ref{eq n3}). 

Next we consider the case $r_0 \leq 2 f(z_0)$.  Since $h$ is positive 
harmonic inside the ball $\B (\vv_0 , (r_0 - f(z_0))/2)$, it is 
a mixture of Poisson kernels which have bounded derivatives inside 
$\B (\vv_0 , (r_0 - f(z_0)) / 4)$.  Thus the maximum of $|\grad h|$
inside $\B (\vv_0 , (r_0 - f(z_0)) / 4)$ is bounded by a constant
times the maximum value of $h$ on $\B (\vv_0 , (r_0 - f(z_0)) / 2)$
divided by the radius of the ball.  From Lemma~\ref{lem h bounds} 
we obtain
$$ |\grad h(\vv_0)| \leq c {\wqt(L) \over \wqt(z_0)} {\log (r_0 / f(z_0)) \over
   \log (g(z_0))} {1 \over r_0 - f(z_0)} \, . $$
Let $b$ denote $r_0 / f(z_0) - 1$.  Then
$$|\grad h (\vv_0)| \leq c {\wqt(L) \over \wqt(z_0)} {\log (1+b) \over \log
   (g(z_0))} {1 \over b f(z_0)} \leq c' {\wqt(L) \over \wqt(z_0)} {1 \over
   r_0 \log (g(z_0))} \, ,$$
finishing the proof in this case.  

It remains to consider the case when both $r_0$ and $z_0$
are at most $2r_f$. Let 
$$\chi_1(r,z) = \max\left (
 {1 \over \wqt(z)} {\log (r / f(z)) \over \log g(z)} ,
{1 \over \wqt(r)}\right)$$
and let $\chi(\xx)$ be the harmonic function
in $\B_{4r_f}\setminus \C$ which is equal to $\chi_1$
on $\partial \B_{4r_f}$ and $0$ on $\partial \C$.
Recall that we have assumed at the beginning
of the section that the boundary of $\C$ is $C^2$-smooth.
It is a standard result that $\chi(\xx)$ is bounded
by a constant (depending on $f$)
times the distance of $\xx$ from $\partial \C$,
for $\xx \in \B_{2r_f}\setminus \C$. By Lemma~\ref{lem h bounds}
we have $h(\xx) \leq \wt c \wqt(L) \chi(\xx) 
\leq c \wqt(L) \hbox{dist}(\xx, \partial \C)$.
We now apply the same argument as in the previous
paragraph to obtain the bound 
$|\grad h(\vv_0)|\leq c' \wqt(L) \hbox{dist}(\xx, \partial \C)
/\hbox{dist}(\xx, \partial \C)= c' \wqt(L)$.
$\Cox$

To prove Lemma~\ref{lem U bound}, we will need the following 
spherical integral.  
\begin{lem} \label{lem sphere}
Let $A_\theta$ be the region on the $s$-sphere defined by 
$$A_\theta = \partial \B_s \setminus (\C \cup \C_{\theta}).$$
Let $r (\xx)$ (respectively $r'(\xx)$) denote the distance between $\xx$
and the $z$-axis (respectively the axis of $\C_{\theta}$).
Then there are constants $\kappa_1 , \kappa_2$ independent of $f$ such 
that for any $s$,
\begin{equation} \label{eq log theta}
\int_{\B_s} {1 \over r (\xx) r'(\xx)} dS 
\leq \kappa_1 + \kappa_2 |\log \theta| ,
\end{equation}
where $dS$ is (non-normalized) area measure on $\partial \B_s$.
Alternatively, 
\begin{equation} \label{eq r}
\int_{A_\theta} {1 \over r (\xx) r'(\xx)} dS \leq 16 \pi \log (\pi g (s))
\end{equation}
independently of $\theta$.  
\end{lem}

\noindent{\sc Proof}: The axis of $\C$ intersects $\partial \B_s$ in two
points, call them $\pp$ and $\wt{\pp}$.  There is an arc in 
$\partial \B_s$ of length $s \theta$ connecting $\pp$ to one of the 
intersection points of $\C_{\theta}$ with $\partial \B_s$; call
this point $\pp'$.  Let $\ww$ denote the midpoint of the arc $\pp \pp'$.
By symmetry through the origin, we may integrate over the
set of points making angle at most $\pi / 2$ with $\ww$, and then double
the result.  Let $\gamma_\uu (\vv)$ denote the arclength along $\partial \B_s$
between the points $\vv$ and $\uu$.  Break the integral 
in~(\ref{eq log theta}) into two pieces:
$$ \int_{\B_s} {1 \over r (\xx) r'(\xx)} dS = 
2 \left [ \int_{\xx : 
   \gamma_\ww (\xx) \leq \theta s} {1 \over r (\xx) r' (\xx) } \, dS +
   \int_{\xx : \theta s \leq \gamma_\ww (\xx) \leq \pi / 2} 
   {1 \over r(\xx) r' (\xx) } \, dS \right ] .$$
When $\gamma_\ww (\xx) \geq \theta s$, then each of $\gamma_\pp (\xx)$ 
and $\gamma_{\pp'} (\xx)$ is at least $\gamma_\ww (\xx) / 2$, and so 
$r (\xx)$ and $r' (\xx)$ are at least $\gamma_\ww (\xx) / 4$.  The integrand
in the second integral is therefore at most $16 / \gamma_\ww (\xx)^2$.  Since
the area of $\{ \xx : a \leq \gamma_\ww (\xx) \leq a + da \}$ is
at most $2 \pi a \, da$, we may integrate over the parameter
$r = \gamma_\ww (\xx)$ to see that the second integral is at most
$$\int_{\theta s}^{\pi s / 2} {32 \pi \over r} \, dr \leq 
   32 \pi (\log (\pi / 2) + |\log \theta|) .$$
To evaluate the first integral, we may integrate over the region
where $\gamma_\pp (\xx) \leq \gamma_{\pp'} (\xx)$ and then double.  On
this region $\gamma_{\pp'} (\xx) \geq \theta s / 2$.  Integrating
over the parameter $r = \gamma_\pp (\xx)$, the first integral is at most
$$ 2 \int_0^{ \theta s } { 2 \over \theta s r} 2 \pi r \, dr
   \leq 8 \pi .$$
Putting these two pieces together proves~(\ref{eq log theta}).

To prove~(\ref{eq r}), use  Cauchy-Schwartz to see that
$$\int_{A_\theta} {1 \over r (\xx) r'(\xx)} dS \leq 
\left( \int_{A_\theta} {1 \over r (\xx)^2} dS \right)^{1/2}
\left( \int_{A_\theta} {1 \over r' (\xx)^2} dS \right)^{1/2}
=\int_{A_\theta} {1 \over r (\xx)^2} dS .  $$
An upper bound for this is 
$$2 \int_{\xx : f(s) / 2 \leq \gamma_\pp (\xx) \leq \pi s / 2}
   {4 \over \gamma_\pp (\xx)^2} \, dS ,$$
which is at most
$$2 \int_{f(s) / 2}^{\pi s / 2} {8 \pi / r} \, dr
    \leq 16 \pi \log (\pi g (s)) .$$
$\Cox$

\noindent{\sc Proof of Lemma}~\ref{lem U bound}: 
We operate by induction on $L$.  First, note that for any $R$ 
and any $L \leq R$, $U(L , \theta) \leq \wqt(L)^{-2}$.
Thus if we choose $b(r) \geq \wqt(4r)^{-2}$, then the result holds for any 
$L \in [R/4 , R]$.  The induction step assumes the result for $L/4$ and proves
the result for $L$.  If $L$ is not regular for $\theta$ then $U(L , \theta) 
\leq U(L/4 , \theta)$ so the induction is trivial.  Thus we may assume $L$ 
is regular.  Applying Lemmas~\ref{lem greens} and~\ref{lem G and q} shows
that for any $L \geq R \geq 2r_f$,
\begin{eqnarray*}
\wqt(L , \theta) & \leq & 2 \left [ \wqt(L)^2 + 10 \int_{\B_{L/4} \setminus 
   (\B_R \cup \C \cup \C_{\theta})}  |\grad h_1| |\grad h_2| 
   \wqt(|\xx|/2 , \theta) |\xx|^{-1} \, d\xx \right. \\
&& + \left. \int_{\B_R \setminus (\C \cup \C_{\theta})} |\grad h_1|
   |\grad h_2| |\xx|^{-1} \, d\xx \right ] .  
\end{eqnarray*}
Write this as an iterated integral, over spherical shells; apply the 
bounds on $|\grad h_i|$ from Lemma~\ref{lem grad bounds} replacing
$z$ by $\rho$ in~(\ref{eq grad bound 1}) at a cost of a factor of at
most some function $\beta (r_f)$, to get
\begin{eqnarray*}
\wqt(L , \theta) & \leq & 2 \beta(r_f)^2 \left [ \wqt(L)^2 + \int_R^{L/4} K_f^2
   {\wqt(L)^2 \over \wqt(s)^2 \log^2 g(s)} {\wqt(s/2 , \theta) \over s}
   \left ( \int_{\B_s \setminus (\C \cup \C_{\theta})}
   {1 \over r(\xx) r'(\xx)} dS \right ) \, ds \right. \\[2ex]
&& + \left. \int_{\sqrt{2} r_f}^R K_f^2 {\wqt(L)^2 \over \wqt(s)^2 \log^2 g(s)} 
   {1 \over s}
   \left ( \int_{\B_s \setminus (\C \cup \C_{\theta})}
   {1 \over r(\xx) r'(\xx)} \, dS \right ) \, ds \right ] \\[2ex]
&& +  \int_{\B_{\sqrt{2}r_f} \setminus (\C \cup \C_{\theta})} |\grad h_1|
   |\grad h_2| |\xx|^{-1} \, d\xx .
\end{eqnarray*}
The last integral is bounded by $\Xi (r_f) \wqt(L)^2$, by 
Lemma~\ref{lem grad bounds}.
Let $R_f$ be large enough so that $\wqt(s) \geq \wqt(s/4) /2$ 
for $s \geq R_f / 4$.  
Change variables in the first line to $t = s/4$ and 
regroup the part where $t < R$ with the second line to get
\begin{eqnarray*}
&& \wqt(L , \theta) \\[2ex]
& \leq & 2 \beta(r_f)^2
   \left [ \wqt(L)^2 + \int_R^{L/16} K_f^2 {\wqt(L)^2
   \over (1/4) \wqt(t)^2 \log^2 g(t)} {\wqt(t , \theta) \over 2t}
   \left ( \int_{\B_s \setminus (\C \cup \C_{\theta})}
   {1 \over r(\xx) r'(\xx)} dS \right ) \, (2 dt) \right. \\[2ex]
&& + \left. \int_{\sqrt{2} r_f}^R K_f^2 {\wqt(L)^2 \over \wqt(s)^2 \log^2 g(s)} 
   {5 \over s}
   \left ( \int_{\B_s \setminus (\C \cup \C_{ \theta})}
   {1 \over r(\xx) r'(\xx)} \, dS \right ) \, ds \right ] \\[2ex]
&& + \, \Xi (r_f) \wqt(L)^2,
\end{eqnarray*}
where the 5 comes from bounding $\wqt(t , \theta)$ above by one, and adding
the regrouped part, which has a total factor of 4.  
Use the first bound from Lemma~\ref{lem sphere} for the
inner integral in the first line and the second bound from the
lemma in the inner integral in the second line and divide by $\wqt(L)^2$
to get
\begin{eqnarray*}
U(L , \theta) & \leq & 8 \beta(r_f)^2 
   (\kappa_1 + \kappa_2 |\log \theta|) K_f^2 
   \int_R^{L/16} {U(t , \theta) \over t \log^2 g(t)} \, dt \\
& + & 2\beta(r_f)^2 + 10 K_f^2 \beta(r_f)^2 
   \int_{\sqrt{2} r_f}^R {16 \pi \log (\pi g (s)) \over 
   s \wqt(s)^2 \log^2 g(s)} ds \\
 & + &  \Xi (r_f).
\end{eqnarray*}
Setting 
$$c_f = 8 \beta (r_f)^2 K_f^2 \max \{ \kappa_1 , \kappa_2 \}$$
and 
$$b(R) = \beta(r_f)^2 \left [ 2 + \Xi (r_f) + 10 K_f^2 \int_1^R {16 \pi 
   \log (\pi g(s)) \over s \wqt(s)^2 \log^2 g(s)} ds \right ] $$
proves the lemma.   $\Cox$

\noindent{\sc Proof of Lemma}~\ref{lem final pieces}: 
The bounds $p_2 < 1$ and $p_1 \leq c_{31} \dd / r_0$ are obvious.

To prove~(\ref{eq n1}), consider three cases.  First suppose that 
$2^k r_0 / 4 \geq L / 24$.  Then 
$z_0 \geq r_0 \geq 2^{-k} L /6$, so
$\wqt(z_0) / \wqt(L) \leq (1+ \ee)^{k+3}$, where $1 + \ee$ is
an upper bound on $\wqt(x) / \wqt(2x)$ for $x \geq r_f$.  Also,
since $z_0 \leq L \leq 6 \cdot 2^k r_0$ and $r_0 / 2 \geq f(z_0)$, 
$${\log g(z_0) \over \log (r_0 / f(z_0))} =
   {\log (z_0 / f(z_0)) \over \log (r_0 / f(z_0))} \leq
   {\log (z_0 / (r_0 / 2)) \over \log (r_0 / (r_0 / 2))} 
= {\log (z_0 / r_0 )+\log 2 \over \log  2}
\leq 2 (k+4) .$$
Thus we obtain
\begin{equation} \label{eq n100}
\sup_{\xx \in A_k} h(\xx) \leq1 \leq  c (k+4) (1 + \ee)^{k+3} 
   {\wqt(L) \over \wqt(z_0)} {\log (r_0 / f(z_0)) \over \log (g(z_0))} .
\end{equation}
Now take an arbitrarily small $\alpha >1$. Then choose
small $\ee>0$ (this requires choosing large $r_f$) and
$\chbounds$ sufficiently large so that 
$c (k+4) (1+\ee)^{k+3}$ is bounded by $\chbounds \alpha^k$.  

The second case is if $z_0 / 2 \leq 2^k r_0 / 4 \leq L / 24$.  This
ensures that $A_k \subseteq \B_{L / 3}$ and thus by 
Lemma~\ref{lem h bounds},
$h(\xx) \leq \wqt(L) / \wqt(|\xx|)$ for any $\xx \in A_k$.  If a point
$\xx \in A_k$ has cylindrical coordinates $z_1$ and $r_1$, then
$$z_1 \leq z_0 + 2^k r_0/4 \leq 2^k r_0/2 + 2^k r_0/4 
\leq 2^k r_0 \leq 2^k z_0$$
and so
$\wqt(z_0) / \wqt(z_1) \leq (1 + \ee)^k$ as in the previous case.
In view of $z_0 \leq 2^k r_0/2$, 
$${\log g(z_0) \over \log (r_0 / f(z_0))} \leq
   {\log (z_0 / r_0 )+\log 2 \over \log  2} \leq
   2 (k+1) .$$
Thus 
$$\sup_{\xx \in A_k} h(\xx) \leq
{\wqt(L) \over \wqt(|\xx|)}(1 + \ee)^k {\wqt(z_1) \over \wqt(z_0)}
{\log (r_0 / f(z_0)) \over \log (g(z_0))}
 2 (k+1) \leq c (k+1) (1 + \ee)^k {\wqt(L) \over \wqt(z_0)}
{\log (r_0 / f(z_0)) \over \log (g(z_0))},$$
which is analogous to~(\ref{eq n100}).

Finally, in the case where $2^k r_0 /4\leq z_0 / 2 \wedge L / 24$, 
let a point $\xx \in A_k$ again have cylindrical coordinates $(z_1 , r_1)$.
Since $z_0 / 2 \leq z_1 \leq 3 z_0 / 2$
and $r_1 \leq 2 \cdot 2^k r_0$, it follows that
$\wqt(z_0) / \wqt(z_1) \leq (1+\ee)^{k+1}$, that $\log g(z_0) / \log g(z_1)
\leq 2$, and that $\log (r_1 / f(z_1)) / \log (r_0 / f(z_0)) \leq
1 + (k+1) \log 2$.  Lemma~\ref{lem h bounds} is again applicable, yielding
$$\sup_{\xx \in A_k} h(\xx) \leq
c {\wqt(L) \over \wqt(z_1)} {\log (r_1 / f(z_1)) \over \log g(z_1)}
(1 + \ee)^k {\wqt(z_1) \over \wqt(z_0)}
{2 \log g(z_1) \over \log g(z_0)}
(1 + (k+1) \log 2)
{\log (r_0 / f(z_0)) \over \log (r_1 / f(z_1))}.$$
This simplifies 
again to~(\ref{eq n100}).  

Recall that $\alpha$ can be chosen arbitrarily close to 1 by choosing 
$\chbounds$ sufficiently large in each of the three cases.  Choosing 
$\alpha < p_2^{-1}$ and $\chbounds$ to be the maximum of the
three values proves~(\ref{eq n1}).

Next we prove~(\ref{eq n4}). Assume that $z_0 \leq r_0$
and find a point $\wt \vv_0$ with the same $\rho$ as for $\vv_0$
and such that $\wt z_0 = \wt r_0$ and $|\vv_0 - \wt \vv_0| < r_0$.
Then $A_k \subset \wt A_{k+4}$ and 
we obtain from~(\ref{eq n1})
$$\sup_{\xx \in  A_{k}} h(\xx) \leq
\sup_{\xx \in \wt A_{k+4}} h(\xx)  \leq  c_{4} \alpha^{k+4} 
{\wqt(L) \over \wqt(\wt z_0)}
   {\log (\wt r_0 / f(\wt z_0)) \over \log (g(\wt z_0))} .$$
Since $\rho/2 \leq \wt z_0 =\wt r_0 \leq r_0 \leq \rho$ we have
$r_0/\rho \geq c$, $\wqt(\wt z_0) \geq  \wqt(\rho)$, 
$\log g(\wt z_0) \geq c\log g(\rho)$, and
$$\log (\wt r_0 / f(\wt z_0)) \leq c\log ( r_0 / f(r_0))
= c\log g(r_0),$$
for some absolute constant $c$. Hence,
$$\sup_{\xx \in  A_{k}} h(\xx) \leq
c' \alpha^k 
{r_0\over \rho} {\wqt(L) \over \wqt(\rho)}
   {\log g (r_0) \over \log (g(\rho))},
$$
which is~(\ref{eq n4}). 

It remains to prove~(\ref{eq n3}) and~(\ref{eq n5}).  
When $r_0 \leq z_0$, scaling down by a factor of 
$2^k r_0$ turns $D_{k+1}$ into a set contained in the union
of two cylinders
with axes at most 1 and radii at most $1 / g(z_0)$, so the
capacity of the rescaled set $D_{k+1}$
is at most a constant multiple of $1 / \log g(z_0)$.
The rescaled point $\yy$ is at distance at least $1/16$ from the
rescaled $D_{k+1}$ and at distance at most 1 from the rescaled $\K$,
so the probability of hitting $D_{k+1}$ before $\K$ starting from
$\yy$ is at most $c / \log g(z_0) \leq c / \log (r_0 / f(z_0))$.
In the case $z_0 \leq r_0$ we use the bound $f(r_0)/r_0 = 1 /g(r_0)$
for the cylinder radius.
$\Cox$

\end{document}